\documentclass[lettersize,journal]{IEEEtran}

\usepackage{cite}
\usepackage{amsmath,amsfonts,amssymb}
\usepackage{algorithm}
\usepackage{algorithmic}
\usepackage{graphicx}
\usepackage{array}
\usepackage{stfloats}
\usepackage{enumerate}
\usepackage{mathrsfs}
\usepackage{textcomp}
\usepackage{epstopdf}  
\usepackage{color}
\usepackage[colorlinks,linkcolor=blue,citecolor=blue]{hyperref} 
\usepackage{multirow}
\usepackage{booktabs}
\usepackage{pifont}
\usepackage{makecell}
\usepackage{caption,subcaption}

\newtheorem{theorem}{Theorem}

\newtheorem{lemma}{Lemma}
\newtheorem{assumption}{Assumption}
\newtheorem{remark}{Remark}
\newtheorem{corollary}{Corollary}

\newcommand{\bsx}{\boldsymbol{x}}
\newcommand{\bse}{\boldsymbol{e}}
\newcommand{\bsv}{\boldsymbol{v}}
\newcommand{\bsu}{\boldsymbol{u}}

\newcommand{\bss}{\boldsymbol{s}}
\newcommand{\bsX}{\boldsymbol{X}}
\newcommand{\bsY}{\boldsymbol{Y}}
\newcommand{\bsy}{\boldsymbol{y}}
\newcommand{\bsz}{\boldsymbol{z}}
\newcommand{\sumn}{\sum\limits_{i=1}^n}
\newcommand{\sumjn}{\sum\limits_{j=1}^n}
\newcommand{\sumT}{\sum\limits_{t=1}^T}

\allowdisplaybreaks[3]

\def\BibTeX{{\rm B\kern-.05em{\sc i\kern-.025em b}\kern-.08em
    T\kern-.1667em\lower.7ex\hbox{E}\kern-.125emX}}
\usepackage{balance}

\begin{document}
\title{ Distributed Online Stochastic Convex-Concave Optimization: Dynamic Regret Analyses under Single and Multiple Consensus Steps}
\author{Wentao Zhang,
        Baoyong~Zhang,~\IEEEmembership{Senior Member,~IEEE},
        Deming~Yuan,~\IEEEmembership{Senior Member,~IEEE},\\
        Shengyuan~Xu,~\IEEEmembership{Senior Member,~IEEE},
        Vincent~K.~N.~Lau,~\IEEEmembership{Fellow,~IEEE}
\thanks{
\emph{Corresponding author: Baoyong Zhang.}}
\thanks{Wentao Zhang is with School of Automation,  Nanjing University of Science and Technology, Nanjing 210094, Jiangsu, P. R. China, and also with the Department of Electronic and Computer Engineering, The Hong Kong University of Science and Technology, Kowloon 999077, Hong Kong (e-mail: iswt.zhang@gmail.com).}
\thanks{Baoyong Zhang, Deming Yuan and Shengyuan Xu are with School of Automation,  Nanjing University of Science and Technology, Nanjing 210094, Jiangsu, P. R. China (e-mail: baoyongzhang@njust.edu.cn, dmyuan1012@gmail.com, syxu@njust.edu.cn).}
\thanks{Vincent K. N. Lau is with the Department of Electronic and Computer Engineering, The Hong Kong University of Science and Technology, Kowloon 999077, Hong Kong (e-mail: eeknlau@ust.hk).
}}
\markboth{IEEE TRANSACTIONS TEMPLATE,~Vol.~XX, JUNE~2025}
{Zhang \MakeLowercase{\textit{et al.}}: Distributed Online Stochastic Convex-Concave Optimization: Dynamic Regret Analyses under Single and Multiple Consensus}
\maketitle
\begin{abstract}
This paper considers the distributed online convex-concave optimization with constraint sets over a multiagent network, in which each agent autonomously generates a series of decision pairs through a designable mechanism to cooperatively minimize the global loss function. To this end,  under no-Euclidean distance metrics, we propose a distributed online stochastic mirror descent convex-concave optimization algorithm with time-varying predictive mappings.
Taking dynamic saddle point regret as a performance metric, it is proved that the proposed algorithm achieves the regret upper-bound in $\mathcal{O}(\max \{T^{\theta_1},  T^{\theta_2} (1+V_T ) \})$ for the general convex-concave loss function,  where $\theta_1, \theta_2 \in(0,1)$ are the tuning parameters,   $T$ is the total iteration time, and $V_T$ is the path-variation. Surely, this algorithm guarantees the sublinear convergence, provided that $V_T$ is sublinear.
 Moreover, aiming to achieve better convergence, we further investigate a variant of this algorithm by employing the multiple consensus technique. The obtained results show that the appropriate setting can effectively tighten the regret bound to a certain extent. Finally, the efficacy of the proposed algorithms is validated and compared through the simulation example of a target tracking problem.
\end{abstract}

\begin{IEEEkeywords}
Distributed  optimization, online convex-concave optimization, Bregman divergence, multiple consensus iterations, dynamic regret.
\end{IEEEkeywords}

\section{Introduction}\label{s1}
\IEEEPARstart{O}{line} convex optimization (OCO) has emerged as a potent methodology that addresses real-time decision-making tasks and has recently attracted extensive attention due to the important applications in smart grids, signal processing, machine learning,
etc \cite{shalev2011online, hazan2016introduction, yuan2024multiO, li2023survey}.  Under the OCO framework, the loss function is only revealed from the adversary after the decision maker commits a decision. By using this function information, the decision maker updates the next decision, thereby generating a series of decisions to achieve the goal of minimizing the cumulative loss function over time. The seminal work on OCO can be traced back to \cite{zinkevich2003online}, where Zinkevich analyzed online gradient descent optimization algorithm and established the regret bound in $\mathcal{O}(\sqrt{T})$. To date, a variety of impressive algorithms have been developed for solving the OCO problem (see, e.g., \cite{cao2019online, 9462561, yuan2020distributed, yi2020distributed, xiong2022distributed, CaoX2022DisCon, NazariP2022DAdam,9239886, shahrampour2017distributed}).

However, the loss functions involved in some important scenarios, such as the bilinear matrix game\cite{beznosikov2020distributed,kovalev2022accelerated}, robust optimization problem \cite{akimoto2021saddle}, transmission and jamming optimization\cite{ChenJ2012convergenceSPP}, constrained optimization duality \cite{yi2020distributed,zhou2019adaptive}, do not apply to the OCO framework but present a convex-concave optimization structure.
Naturally, these practical scenarios spark an interest in exploring online convex-concave optimization (OCCO), also known as online saddle point problems.
In this paper, we investigate a distributed solution of OCCO over a multiagent network, that is, solving the specific optimization problem formulated in (\ref{problem definition}).
\begin{flalign}\label{problem definition}
\min\limits_{\boldsymbol{x}_t \in \boldsymbol{X}} \max\limits_{\boldsymbol{y}_t \in \boldsymbol{Y}}\,\,\sum\limits_{t=1}^T f_t ( \boldsymbol{x}_t,\boldsymbol{y}_t )
\end{flalign}
where  $f_t( \boldsymbol{x}_t,\boldsymbol{y}_t )= \sum_{i=1}^nf_{i,t} ( \boldsymbol{x}_t,\boldsymbol{y}_t )$, $f_{i,t}(\cdot, \cdot): \mathbb{R}^d \times \mathbb{R}^m\!\rightarrow\! \mathbb{R}$ are the local convex-concave functions known only to agent $i$ at time $t$, and $\boldsymbol{X} \subset \mathbb{R}^{d}, \boldsymbol{Y} \subset \mathbb{R}^{m}$ are two convex and compact sets with
$
 \max_{\boldsymbol{x},\boldsymbol{x}_a \in \boldsymbol{X}} \| \boldsymbol{x} -\boldsymbol{x}_a \| \leq M_X ,
 \max_{\boldsymbol{y},\boldsymbol{y}_a \in \boldsymbol{Y}} \| \boldsymbol{y} -\boldsymbol{y}_a \| \leq M_Y, M_X>0, M_Y>0$.

%
%

In the gradient feedback of OCCO at each round, we consider the stochastic gradient from practical application scenarios, which simultaneously introduces randomness. Therefore, we utilize the  \emph{expected dynamic saddle point regret} in (\ref{SP-Regret-j}) as a performance metric, that is,
\begin{flalign} \label{SP-Regret-j}
&\textbf{ESP-Regret}_{d}^j (T)  =  \left|  \mathbb{E} \left[ \sum\limits_{t=1}^T f_t (\boldsymbol{x}_{j,t},\boldsymbol{y}_{j,t}) -  \sum\limits_{t=1}^T  f_t (\boldsymbol{x}_t^*,\boldsymbol{y}_t^*) \right]\right|
\end{flalign}
to measure the efficiency of the developed algorithm, where the saddle point
$
(\boldsymbol{x}_t^*,\boldsymbol{y}_t^*) \in {\arg\min}_{\boldsymbol{x} \in \boldsymbol{X}}    {\arg\max}_{\boldsymbol{y} \in \boldsymbol{Y}}  f_t (\boldsymbol{x},\boldsymbol{y})
$
 of problem (\ref{problem definition}) satisfies the property
$f_t ( \boldsymbol{x}_t^*,\boldsymbol{y} )  \leq f_t ( \boldsymbol{x}_t^*,\boldsymbol{y}_t^* ) \leq f_t ( \boldsymbol{x},\boldsymbol{y}_t^* ), t\in [T].$

 \begin{table*}[t] 
 \label{table_1}
 \renewcommand\arraystretch{1.8}
\begin{center}
  \caption{The comparison of related researches on OCCO.}
  \resizebox{\linewidth}{!}{ 
  \label{table_1}
 \begin{tabular}{ccccccc}
 \toprule
 References  & Loss function & \makecell[c]{Distributed\\ manner}& \makecell[c]{Non-Euclidean\\ space} &\makecell[c]{Stochastic \\ gradient} &\makecell[c]{ Performance\\Metric} & Regret bound over $T^{\dag}$  \\ \midrule
Ho-Nguyen \emph{et al.}\cite{ho2019exploiting}  &\makecell[c]{ Convex-concave\\ Non-smooth} &\ding{55}  &\ding{51} &\ding{55}  & Weighted gap & $\mathcal{O}(\sqrt{T})$  \\
Rivera \emph{et al.} \cite{rivera2018online}  & \makecell[c]{ Convex-concave\\ Non-smooth}  & \ding{55}  &\ding{55} &\ding{55}  & Static regret  & $\mathcal{O}(\sqrt{T} \ln T)$ \\
%
 Wood \emph{et al.} \cite{wood2023online}   & SC-SC; Smooth &\ding{55} &\ding{55} & \ding{51}& Equilibrium points  & $\mathcal{O}(1+V_T^I)^{\dag}$\\

 Cardoso \emph{et al.}\cite{cardoso2019competing}  & \makecell[c]{Convex-concave \\ Non-smooth }  & \ding{55}  & \ding{55} & \ding{55}  & Static regret & $\mathcal{O}(kT^{5/6}\ln T)$ \\
 Roy \emph{et al.} \cite{roy2019online}   & SC-SC; Smooth & \ding{55}    &\ding{55} & \ding{51} & Dynamic regret  &  $\mathcal{O}(\sqrt{T \max \{W_T,M_T \}} )$  \\[3pt]
Zhang \emph{et al.} \cite{10239326}  & \makecell[c]{Convex-concave \\ Non-smooth } & \ding{51}  &\ding{55}   &\ding{55}   & Dynamic regret & \makecell[c]{$\mathcal{O}\left(\max \{T^{a_1},  T^{a_2} (1+V_T^I ) \}\right)$ \\[2pt]$a_1=\max \{1-b_1, 1-b_2 \}, a_2=\max \{b_1, b_2 \},b_1,b_2\in(0,1)$} \\[6pt]
 \textbf{Algorithm \ref{algorithm 1}}   & \makecell[c]{Convex-concave \\ Non-smooth } & \ding{51} & \ding{51} & \ding{51} & Dynamic regret  & \makecell[c]{$\mathcal{O}\left(\max \{T^{\theta_1},  T^{\theta_2} (1+V_T ) \}\right)$ \\[2pt]$\theta_1=\max \{1-\gamma_1, 1-\gamma_2 \}, \theta_2=\max \{\gamma_1, \gamma_2 \},\gamma_1,\gamma_2\in(0,1)$}  \\[6pt]
\textbf{Algorithm \ref{algorithm 2}}   & \makecell[c]{Convex-concave \\ Non-smooth } & \ding{51} & \ding{51} & \ding{51} & Dynamic regret  & \makecell[c]{$\mathcal{O}\left( \max\left\{ \left(1+\frac{\Gamma_1{\sigma_1}^{\underline{K}-1}}{1-{\sigma_1}^{\underline{K}}}\right)T^{\theta_1},  T^{\theta_2} (1+V_T ) \right\} \right)$ \\[4pt] $\sigma_1\in(0,1), \underline{K}=\min_{t\in[T]} \{ K_t\}, K_t \in \mathbb{Z}_+$}  \\[6pt] \bottomrule
\end{tabular}
}
\begin{flushleft}\scriptsize
$\dag$ Note: $W_T$ and $M_T$ are a function variation and path variation defined in \cite{roy2019online}, respectively. $V_T^I$ represents the variation $V_T$ satisfying $B_t=I_d$ and $C_t=I_m$.
\end{flushleft}
\end{center}
\vspace{-2.5em}
\end{table*}

Commonly, the analysis of dynamic regrets depends on specific features of the online optimization problem\cite{besbes2015non}. In light of this, incorporating time-varying prediction mappings $B_t$ and $C_t$, we utilize path-variations defined in (\ref{assumption Path Variation equation}) to elucidate the variation degree between optima.
\begin{flalign} \label{assumption Path Variation equation}
V_T^x:= \sum_{t=1}^T \| \boldsymbol{x}_{t+1}^* -B_t \boldsymbol{x}_t^* \|, \
 V_T^y:= \sum_{t=1}^T \| \boldsymbol{y}_{t+1}^* - C_t\boldsymbol{y}_t^* \|.
\end{flalign}
Note that (\ref{assumption Path Variation equation}) naturally covers the regular path-variations, i.e., the case with $B_t=I_d, C_t=I_m$ (see \cite{10239326,9239886,zhang2019distributed}).
Moreover, in certain application scenarios characterized by a dynamic relationship between optima, like the target tracking problem discussed in \cite{shahrampour2017distributed}, the efficient predictive mappings have the capability to establish the small $V_T^x$ and $V_T^y$.
Thus, compared with the regular form, (\ref{assumption Path Variation equation}) is more general.
To facilitate the following analysis, denote $V_T =\max \{V_T^x, V_T^y\}$.

The objective of this paper is to design an effective distributed online convex-concave algorithm such that the dynamic saddle point regret (\ref{SP-Regret-j}) grows sublinearly.
\subsection{Literature Review}
\allowdisplaybreaks{The research in this paper is related to two bodies of literature: centralized solutions ($n=1$) and distributed solutions for OCCO. The overview of the related works is stated below.

In \cite{ho2019exploiting}, Ho-Nguyen and K{\i}l{\i}n\c{c}-Karzan earlier investigated the centralized OCCO and established the sublinear convergence for their proposed algorithm in a metric of weighted online saddle point gap.
The work \cite{rivera2018online} studied an algorithm named online saddle point follow-the-leader, and for the general convex-concave loss function, it attained  the  static regret in order $\mathcal{O}(\sqrt{T} \ln T)$. Subsequently, Xu \emph{et al.} \cite{xu2019online} additionally considered a regularization term based on \cite{rivera2018online} to enhance the decision quality. The work \cite{wood2023online} conducted an analysis into a class of OCCO with decision-dependent data, employing the theory of equilibrium points. Focusing on the specific scenarios of OCCO,  Cardoso \emph{et al.}\cite{cardoso2019competing} investigated the Nash Equilibrium regrets of online zero-sum game with full-information and bandit feedbacks.
In \cite{roy2019online}, Roy \emph{et al.} developed online extragradient and Frank-Wolfe convex-concave optimization (CCO) algorithms and showed that their sublinear convergence under two regret metrics. However, it is essential to acknowledge that these results rely on the stringent assumptions that the loss function possesses strongly convex-strongly concave (SC-SC) characteristics and is smooth.

It is well known that centralized algorithms may be limited and powerless for large-scale optimization problems and complex scenarios due to the computational bottleneck of a processor. In contrast, the distributed algorithms over a multiagent network overcome this limitation and have attracted the attention of many researchers (see, e.g., \cite{yang2019survey,nedic2018distributed,li2023survey,nedic2008distributed,LiuC2024distributed,NiuD2025adual, YangZ2025Differentially, xu2021distributed, zhang2020distributed, yi2020distributed,9239886}).
For distributed off-line CCO, the work\cite{mateos2016distributed} proposed a projected subgradient algorithm with Laplacian averaging for the cases with explicit agreement constraints and analyzed its applications on distributed convex optimizations.  The work\cite{rogozin2021decentralized} showed the lower bounds for distributed smooth CCO and studied distributed mirror-prox convex-concave algorithms. With similar ideas, the work \cite{beznosikov2020distributed} studied the lower bounds under the stochastic condition. Considering similarities between local loss functions, Beznosikov \emph{et al.} \cite{beznosikov2021distributed} investigated the min-max data similarity algorithms under centralized and distributed networks and obtained their communication complexity bounds.
Qureshi \emph{et al.}\cite{qureshi2023distributed} studied a distributed stochastic gradient method with gradient tracking and showed its linear convergence with an error neighborhood for strongly concave-convex functions.
However, the above mentioned distributed off-line algorithms are difficult to handle OCCO because the loss function is time-varying and unknowable in advance \cite{hazan2016introduction, besbes2015non,yuan2020distributed}.
For this case, Zhang \emph{et al.} \cite{10239326} designed two distributed online subgradient saddle point optimization algorithms under two information feedbacks, and the related results showcased the effectiveness of these algorithms by obtaining sublinear dynamic saddle point regrets.
A detailed comparison about OCCO is shown in Table \ref{table_1}.

Note that in terms of the distributed solutions of OCCO, the existing research is far from sufficient. In addition, the distributed algorithms developed in \cite{10239326} are difficult to fully exploit certain properties depending on the optimization problem due to the use of the traditional Euclidean distance. For example, under Kullback-Leibler (KL) divergence, the mathematical equations for the solutions to the optimization problem with simplex constraint are explicitly available, whereas under Euclidean distance, these are inaccessible and the solutions depend on solving projection operations \cite{yuan2020distributed}.  Moreover, considering that true gradients are often difficult to obtain due to measurement errors and inaccurate calculations \cite{polyak1987introduction}, we use stochastic gradients instead of it in the proposed algorithm, which results in a more general and practical version.

On the other hand, the frequency of information exchange among agents at each time $t$ significantly effects the consensus process of distributed algorithms. In contrast to single consensus, the distributed algorithms with multiple consensus at time $t$ can enhance information diffusion across the network, which enables each agent to aggregate localized decisions from agents further away to accelerate global consensus and improve optimization efficiency \cite{eshraghi2020distributed,jakovetic2014fast}. In addition, for limited communication networks, such as network delays and noise, multiple consensus has better robustness than single one. Based on the above points, the distributed stochastic algorithm with  single and multiple consensus for OCCO under a non-Euclidean distance is well-motivated.}


\vspace{-0.25cm}
\subsection{Contributions}
The contributions of this paper are summarized as follows.

 1) Taking Bregman divergence as a generalized distance metric, we propose a distributed solution with the optional predictive mappings in a non-Euclidean sense for OCCO. Benefiting from the free selectivity of Bregman divergence, this solution is more flexible for different optimization problems than the one with the traditional Euclidean distance in \cite{10239326}. In addition, the use of the predictive technique can further improve the quality of the committed decisions at each round based on the factitious knowledge and experience.

 2) By combining the mirror descent method and time-varying predictive technique, a \emph{distributed online stochastic mirror descent convex-concave optimization }(DOSMD-CCO) algorithm is developed, in which stochastic gradients are used to addressing the inaccuracy in obtaining true gradients.
 For the convex-concave loss function, we show its expected dynamic saddle point regret scaling in $\mathcal{O}(\max \{T^{\theta_1},  T^{\theta_2} (1+V_T ) \})$, where $\theta_1, \theta_2 \in(0,1)$ are two tuning scalars. Clearly, the developed algorithm can guarantee the sublinear dynamic regret with respect to $T$ under the premise of sublinear $V_T$ and allow a potential performance improvement through finding appropriate mappings $B_t$ and $C_t$.

 3) Further, aiming to achieve better convergence performance, we investigate a multiple consensus version of Algorithm DOSMD-CCO by employing the technique of multiple consensus iterations. The theoretical results give the effect of consensus parameters on regret bound and show that this technique can effectively tighten this regret bound to a certain extent. Finally, the effectiveness of the proposed algorithms is validated and comparatively through a simulation example.
\vspace{-0.25cm}
\subsection{Notations}

 ${\mathbb{R}}^n$ and $\mathbb{Z}$ $(\mathbb{Z}_+)$ represent the sets of the $n$-dimensional vectors and (positive) integers, respectively. The symbol $\|\boldsymbol{u}\|_2$ $(\|\boldsymbol{u}\|_1)$  stands for  the Euclidean  norm (1-norm) of a vector $\boldsymbol{u}$. Denote $\|\cdot\|_*$ as the dual norm of $\| \cdot \|$. Denote $I_d$ as $d\times d$ identity matrix.  Let $[A_t]_{ij}$ stand for  the $(i,j)$-th element of matrix $A_t$.
 Write $[\bsx]_i$ and $[m]$ to denote the $i$th entry of vector $\bsx$ and the integer set $\{1,2,\ldots,m \}$, respectively. The simplex $\{\bsx \in \mathbb{R}^n | \sum_{i=1}^n [\bsx]_i=1,[x]_i \geq 0,i \in [n]\}$ is denoted as $\triangle_n$.
 Let $\mathcal{F}_t$ be $\sigma$-field consisting from the entire history information of random variables up to time $t$.

\section{Preliminaries}\label{s2}
In this section, we introduce some preliminaries about $\mathcal{G}_t$, $f_{i,t}$, and Bregman divergence.
\vspace{-0.2cm}
\subsection{Graph Theory and Basic Assumptions}
 Denote $\mathcal{G}_t:=\{\mathcal{V},\mathcal{E}_t,A_t \}$ as a directed time-varying network (graph), in which $\mathcal{V} : = [n] $ and $ \mathcal{E}_t \subseteq \mathcal{V} \times \mathcal{V} $ are the node and edge sets, respectively, and $A_t \in \mathbb{R}^{n \times n}$ is the weighted matrix satisfying doubly stochasticity.  Let $\mathcal{N}_{i}^{\text {in }}(t)=\{i\} \cup\{j \mid(j, i) \in \mathcal{E}_t\}$ represents the in-neighbors of agent $i$.  The weighted matrix fulfills that $[A_{t}]_{ij}>\zeta, 0<\zeta <1$ if $j \in \mathcal{N}_{i}^{\text {in }}(t)$, and  $[A_{t}]_{ij}=0$ otherwise. Around the graph $\mathcal{G}_t$, we firstly give a standard assumption and a basic lemma.
\begin{assumption} \label{graph assumption 1} \cite{nedic2008distributed,yi2020distributed}
There exists a positive integer $Q$ such that for all non-negative integer $k$, the graph $(\mathcal{V}, \bigcup_{i=k Q+1}^{(k+1)Q} \mathcal{E}_{i})$ is strongly connected.
\end{assumption}

\begin{lemma} \label{lemma-graph}\cite{nedic2008distributed} 
 Suppose Assumption \ref{graph assumption 1} hold. Then, we have that for all $i, j \in \mathcal{V}$ and all $ t \geq s \geq 1,$
\begin{flalign}
\left|[\Phi(t, s)]_{i j}-\frac{1}{n}\right| \leq \Gamma \sigma^{(t-s)}
\end{flalign}
where  $\Phi(t, s)=A_{t} A_{t-1} \ldots A_{s}$, $\Gamma=(1-\zeta / 4 n^{2})^{(1-2Q) /  Q}$ and $\sigma=(1-\zeta / 4 n^{2})^{1 /  Q}$.
\end{lemma}


Instead of true gradients, stochastic gradients are considered into the algorithm design, which is more practical and general due to the inaccuracy and measurement errors in obtaining a true gradient \cite{polyak1987introduction}. Suppose that there exist two independent stochastic oracles that can generate the noisy gradients satisfying the conditions in Assumption \ref{assu schostic gradient}.
\begin{assumption}\label{assu schostic gradient}
 \cite{nedic2014stochastic, xiong2022event}  The stochastic gradients of $f_{i,t}$ satisfy that for any $\bsx \in \bsX, \bsy \in \bsY$,
\begin{flalign*}
 (i)\ & \mathbb{E} \left[\widetilde{\nabla}_{i,t}^x (\bsx, \bsy) | \mathcal{F}_{t-1}  \right]= \nabla_x f_{i,t} (\bsx,\bsy), \quad \quad \quad \quad \quad \quad\quad \quad\ \\
 &\mathbb{E} \left[\widetilde{\nabla}_{i,t}^y(\bsx, \bsy) | \mathcal{F}_{t-1}  \right]= \nabla_y f_{i,t} (\bsx,\bsy); \\
 (ii)\ & \mathbb{E} \left[ \|\widetilde{\nabla}_{i,t}^x (\bsx, \bsy)\|_*^2 |  \mathcal{F}_{t-1}  \right] \leq L_X^2,\\
 &\mathbb{E} \left[ \|\widetilde{\nabla}_{i,t}^y(\bsx, \bsy) \|_*^2 | \mathcal{F}_{t-1} \right] \leq L_Y^2.
\end{flalign*}
Denote $\widetilde{\nabla}_{i,t}^x := \widetilde{\nabla}_{i,t}^x (\bsx_{i,t}, \bsy_{i,t}), \widetilde{\nabla}_{i,t}^y := \widetilde{\nabla}_{i,t}^y (\bsx_{i,t}, \bsy_{i,t})$.
\end{assumption}

From Assumption \ref{assu schostic gradient} and law of total expectation, it can be obtained that $\|\nabla_x f_{i,t} (\bsx,\bsy)\|_* =\|\mathbb{E} [\widetilde{\nabla}_{i,t}^x (\bsx, \bsy)]\|_* \leq\mathbb{E} [\|\widetilde{\nabla}_{i,t}^x (\bsx, \bsy)\|_*] \leq  \{ \mathbb{E} [ \|\widetilde{\nabla}_{i,t}^{x}(\bsx,\bsy) \|_*^2 ]\}^{\frac{1}{2}} \leq L_X.$ Similarly, $\|\nabla_y f_{i,t} (\bsx,\bsy)\|_* \leq L_Y.$  Based on lemma 2.6 in \cite{shalev2011online}, this further implies that for any $\bsx_a, \bsx_b \in \bsX$ and $\bsy_a, \bsy_b \in \bsY$,
\begin{flalign}\label{fun lipschitz}
|f_{i,t}(\boldsymbol{x}_{a}, \boldsymbol{y}_a)& -f_{i,t}(\boldsymbol{x}_{b}, \boldsymbol{y}_b)|\nonumber \\
&\leq L_{X} \|\boldsymbol{x}_{a}-\boldsymbol{x}_{b}\|+L_{Y} \|\boldsymbol{y}_{a}-\boldsymbol{y}_{b}\|.
\end{flalign}

\subsection{Bregman Divergence}
This paper focuses on the development of an online CCO algorithm using the mirror descent approach, in which the Bregman divergence defined in (\ref{definition DR x DR y}) is a central component in this  approach, serving as a distance-measuring function.
\begin{flalign}\label{definition DR x DR y}
&\Psi_{\mathcal{R}}^x (\bsx_a,\bsx_b):=\mathcal{R}^x(\bsx_a)-\mathcal{R}^x(\bsx_b)- \langle \nabla\mathcal{R}^x(\bsx_b), \bsx_a-\bsx_b \rangle, \nonumber\\
&\Psi_{\mathcal{R}}^y (\bsy_a,\bsy_b):=\mathcal{R}^y(\bsy_a)-\mathcal{R}^y(\bsy_b)- \langle \nabla\mathcal{R}^y(\bsy_b), \bsy_a-\bsy_b \rangle
\end{flalign}
where $\mathcal{R}^x: \mathbb{R}^d \rightarrow \mathbb{R}$ and $\mathcal{R}^y: \mathbb{R}^m \rightarrow \mathbb{R}$ are the associated distance-generating functions
and satisfy $\varrho_x$- and $\varrho_y$-strong convexity on the sets $\bsX$ and $\bsY$, respectively. Further, based on the strong convexity, (\ref{definition DR x DR y}) follows
\begin{flalign}\label{DR x DR y lower bound}
\Psi_{\mathcal{R}}^x (\bsx_a,\bsx_b)\geq \frac{\varrho_x}{2} \|\bsx_a-\bsx_b \|^2, \forall \bsx_a,\bsx_b\in\bsX \nonumber\\
\Psi_{\mathcal{R}}^y (\bsy_a,\bsy_b)\geq \frac{\varrho_y}{2} \|\bsy_a-\bsy_b \|^2, \forall \bsy_a,\bsy_b \in \bsY.
\end{flalign}

By selecting diverse distance-generating functions, Bregman divergences can be derived, including notable examples such as the standard Euclidean distance and Kullback-Leibler (KL) divergence (see more examples in \cite{banerjee2005clustering,bauschke2001joint,dhillon2008matrix}).
Considering $\mathcal{R}^x$ and $\mathcal{R}^y$ along with the associated $\Psi_{\mathcal{R}}^x$ and $\Psi_{\mathcal{R}}^y$, we introduce the following standard assumptions for the dynamic regret analysis of online mirror descent algorithms \cite{shahrampour2017distributed, yi2020distributed2,li2021distributed}.

\begin{assumption} \label{assu DR continous}
 i) For all $ \bsx_a, \bsx_b, \bsx_c \in \bsX,$ $ \bsy_a, \bsy_b, \bsy_c \in \bsY$,
\begin{flalign*}
&|  \Psi_{\mathcal{R}}^x(\bsx_a, \bsx_c)-\Psi_{\mathcal{R}}^x(\bsx_b, \bsx_c)| \leq K_{X} \|\bsx_a-\bsx_b \|, \nonumber \\ &|\Psi_{\mathcal{R}}^y(\bsy_a,\bsy_c)-\Psi_{\mathcal{R}}^y(\bsy_b,\bsy_c)| \leq K_{Y} \|\bsy_a-\bsy_b \|.
\end{flalign*}
ii) For $\bsx, \bsz_i \in \mathbb{R}^d$, $\bsy, \bsv_i \in \mathbb{R}^m$,
\begin{flalign*}
&\Psi_{\mathcal{R}}^x \left(\bsx, \sum_{i=1}^n [\boldsymbol{s}_1]_i \boldsymbol{z}_i \right) \leq \sum_{i=1}^n [\boldsymbol{s}_1]_i \Psi_{\mathcal{R}}^x (\bsx,  \boldsymbol{z}_i ), \forall \boldsymbol{s}_1 \in \triangle^n,\nonumber \\
&\Psi_{\mathcal{R}}^y \left(\bsy, \sum_{i=1}^n [\boldsymbol{s}_2]_i \boldsymbol{v}_i \right) \leq \sum_{i=1}^n [\boldsymbol{s}_2]_i \Psi_{\mathcal{R}}^y (\bsy,  \boldsymbol{v}_i ), \forall \boldsymbol{s}_2 \in \triangle^n.
\end{flalign*}
\end{assumption}

\begin{assumption} \label{assu DRX Y noexpansive}
{\setlength{\spaceskip}{0.3em} For all $ \bsx_a,  \bsx_b  \in \bsX, \bsy_a,  \bsy_b\in\bsY,$}
\begin{flalign*}
&\Psi_{\mathcal{R}}^x (B_t \bsx_a, B_t \bsx_b ) \leq  \Psi_{\mathcal{R}}^x (\bsx_a,  \bsx_b ),\nonumber \\
&\Psi_{\mathcal{R}}^y ( C_t\bsy_a, C_t \bsy_b ) \leq  \Psi_{\mathcal{R}}^y (\bsy_a,  \bsy_b )
\end{flalign*}
hold, $\| B_t\| \leq 1, \|C_t \|\leq 1, \forall t \in {T}$, and $B_t\bsx \in \bsX, C_t \bsy \in \bsY$ hold as long as $\bsx\in \bsX, \bsy \in \bsY$.
\end{assumption}

Assumption \ref{assu DRX Y noexpansive} ensures that both mappings $B_t$ and $C_t, t\in[T]$ are nonexpansive and non-violating, i.e., as Algorithm \ref{algorithm 1} progresses, the negative impact of an inaccurate prediction at a specific time does not continuously intensify. The identity mappings satisfy it and a similar requirement also is made in \cite{shahrampour2017distributed, yi2020distributed2}.
\vspace{-0.3cm}
\section{Distributed Online Stochastic Mirror Descent Convex-Concave Optimization Algorithm} \label{s3}

\subsection{ Algorithm Design}
The DOSMD-CCO algorithm  is presented in Algorithm \ref{algorithm 1}, which involves the following key steps.

1) \emph{Mirror descent step in Step 3:} Considering the inaccuracy and measurement errors in obtaining a true gradient, the stochastic gradients are utilized as a practical and general solution.  Based on this, agent $i\in\mathcal{V}$ executes the mirror descent steps to obtain the auxiliary variables $\nabla \mathcal{R}^x (\bsz_{i,t}), \nabla \mathcal{R}^y (\bsv_{i,t})$.

2) \emph{Bregman projections in Step 4:} To guarantee the effectiveness of decision-making, the Bregman projections for variables $\tilde{\bsx}_{i,t}$ and $\tilde{\bsy}_{i,t}$ are executed, respectively, in which the Bregman divergences $\Psi_{\mathcal{R}}^x$ and $\Psi_{\mathcal{R}}^y$ are employed as more flexible distance-measuring functions.

3) \emph{Predictions in Step 5:} Through selecting the appropriate mappings $B_t$ and $C_t$ relying on factitious experience, agent $i$ establishes the corrected decisions $\bss_{i,t}^x$ and $\bss_{i,t}^y$ with predictions, which can be better than the one without predictions.

4) \emph{Consensus in Step 6:} By communicating with its neighbors, agent $i$ receives the auxiliary decisions $\bss_{j,t}^x, \bss_{j,t}^y, j \in \mathcal{N}_{i}^{\text {in }}(t)$ and then executes the consensus steps to output the decisions $\bsx_{i,t+1}$ and $\bsy_{i,t+1}$ at time $t+1$.

\begin{algorithm}[!t]
	\renewcommand{\algorithmicrequire}{\textbf{Initialize:} }
	\caption{  DOSMD-CCO algorithm.}
	\label{algorithm 1}
	\begin{algorithmic}[1]
		\REQUIRE Initial  decisions $\boldsymbol{x}_{i,1} \in \boldsymbol{X} ,\boldsymbol{y}_{i,1} \in \boldsymbol{Y},$ the parameters $\alpha_t,\eta_t>0$, and the mappings $B_t, C_t.$
\FOR {$t=1,2,\cdots,T$}
             \FOR {$i \in \mathcal{V}$ in parallel}

             \STATE Agent $i$ gets the stochastic gradients $\widetilde{\nabla}_{i,t}^x, \widetilde{\nabla}_{i,t}^y$, and computes, respectively,
              \setlength{\parskip}{0.2em}

%

\begin{center}
$\nabla \mathcal{R}^x (\bsz_{i,t})= \nabla \mathcal{R}^x (\bsx_{i,t}) -\alpha_t \widetilde{\nabla}_{i,t}^x,$

$\nabla \mathcal{R}^y (\bsv_{i,t})= \nabla \mathcal{R}^y (\bsy_{i,t}) +\eta_t \widetilde{\nabla}_{i,t}^y.$
\end{center}

\STATE Executes the Bregman projections, respectively,
\begin{center}
$\tilde{\bsx}_{i,t}= \underset{\bsx\in\bsX}{\arg\min} \ \Psi_{\mathcal{R}}^x (\bsx, \bsz_{i,t}),$

$\tilde{\bsy}_{i,t}=\underset{\bsy\in\bsY}{\arg\min}\ \Psi_{\mathcal{R}}^y (\bsy, \bsv_{i,t}).$
\end{center}

\STATE Runs the decision corrections using prediction mappings, i.e.,

\begin{center}
$\bss_{i,t}^x= B_t \tilde{\bsx}_{i,t},\ \bss_{i,t}^y= C_t \tilde{\bsy}_{i,t}.$		
\end{center}
		
\STATE Receives the predictions $\boldsymbol{s}_{j,t}^x$ and $\boldsymbol{s}_{j,t}^y$ from its in-neighbors, and updates decisions by executing 

             \begin{center}
             $\boldsymbol{x}_{i,t+1}= \sum\limits_{j \in \mathcal{N}_{i}^{in}(t)}{[A_t]_{ij} \bss_{j,t}^x}$,

           		$\boldsymbol{y}_{i,t+1}= \sum\limits_{j \in \mathcal{N}_{i}^{in}(t)}{[A_t]_{ij} \bss_{j,t}^y}$.
             \end{center}
		\ENDFOR
		\ENDFOR
	\end{algorithmic}
\end{algorithm}
\vspace{-0.3cm}
\subsection{Main Convergence Results}
We firstly show some necessary lemmas for the regret analysis. Denote $\boldsymbol{g}_{i,t}^{x}:=\nabla_x f_{i,t} (\bsx_{i,t},\bsy_{i,t})$, $ \boldsymbol{g}_{i,t}^{y}:=\nabla_y f_{i,t} (\bsx_{i,t},\bsy_{i,t})$. Write the running averages of $\boldsymbol{x}_{i,t}$, $\boldsymbol{y}_{i,t}$ as
$
\boldsymbol{x}_{avg,t}= \frac{1}{n} \sum_{i=1}^n \boldsymbol{x}_{i,t}, \ \boldsymbol{y}_{avg,t}= \frac{1}{n} \sum_{i=1}^n \boldsymbol{y}_{i,t} , \ \forall  t \in [T].
$

\begin{lemma} \label{lemma xx yy}
{\setlength{\spaceskip}{0.25em} For any $i\in \mathcal{V}$, we have that}
\begin{flalign}
(i)\ &\| \tilde{\bsx}_{i,t}  -{\bsx}_{i,t}\| \leq \frac{\alpha_t } {\varrho_x} \left\|  \widetilde{\nabla}_{i,t}^x \right\|_*,\nonumber \\
(ii)\ &\| \tilde{\bsy}_{i,t}  -{\bsy}_{i,t}\|    \leq \frac{\eta_t}{\varrho_y} \left\|  \widetilde{\nabla}_{i,t}^y \right\|_*.
\end{flalign}
\end{lemma}
\noindent{\em Proof:}
$(i)$ According to the optimality of $\tilde{\bsx}_{i,t}$,  we yield by using the first-order optimality condition that
$\langle \nabla \Psi_{\mathcal{R}}^x (\tilde{\bsx}_{i,t}, \bsz_{i,t}), \bsx- \tilde{\bsx}_{i,t} \rangle \geq0,\  \forall \bsx \in \bsX.$
Due to $\nabla \Psi_{\mathcal{R}}^x (\bsx, \bsz_{i,t})= \nabla \mathcal{R}^x (\bsx) - \nabla \mathcal{R}^x(\bsz_{i,t})$, it follows
\begin{flalign} \label{lem xx yy eq01}
&\langle \nabla \mathcal{R}^x ( \tilde{\bsx}_{i,t} ) - \nabla \mathcal{R}^x(\bsz_{i,t}), \bsx- \tilde{\bsx}_{i,t} \rangle \geq0,\  \forall \bsx \in \bsX.
\end{flalign}

Substituting $\nabla \mathcal{R}^x(\bsz_{i,t})$ from Algorithm \ref{algorithm 1} into (\ref{lem xx yy eq01}) and rearranging it, we have that
\allowdisplaybreaks{\begin{flalign}\label{lem xx yy eq02}
&\langle \widetilde{\nabla}_{i,t}^x, \tilde{\bsx}_{i,t} -\bsx \rangle \nonumber \\
&\leq \frac{1}{\alpha_t}\left\langle \nabla \mathcal{R}^x ( \tilde{\bsx}_{i,t} ) - \nabla \mathcal{R}^x(\bsx_{i,t}), \bsx- \tilde{\bsx}_{i,t} \right\rangle \nonumber \\
&\overset{(a)}{=} \frac{1}{\alpha_t}\left[ \Psi_{\mathcal{R}}^x(\bsx, \bsx_{i,t}) -\Psi_{\mathcal{R}}^x(\bsx, \tilde{\bsx}_{i,t})- \Psi_{\mathcal{R}}^x(\tilde{\bsx}_{i,t}, \bsx_{i,t}) \right]
\end{flalign}
where $(a)$ is obtained by using the fact that $\langle \nabla R^x (\bsx)- \nabla R^x (\bsz), \bsy- \bsz \rangle= \Psi_{\mathcal{R}}^x(\bsy, \bsz)+ \Psi_{\mathcal{R}}^x(\bsz, \bsx)-\Psi_{\mathcal{R}}^x(\bsy, \bsx)$.}

 Further, let $\bsx=\bsx_{i,t} \in \bsX$, and we can obtain by using the facts $\alpha_t>0$ and (\ref{DR x DR y lower bound}) that
\begin{flalign}
\varrho_x \|\bsx_{i,t}- \tilde{\bsx}_{i,t} \|^2    &\leq \langle\alpha_t  \widetilde{\nabla}_{i,t}^x, \bsx_{i,t} -\tilde{\bsx}_{i,t}  \rangle \nonumber\\
&\leq \alpha_t \|  \widetilde{\nabla}_{i,t}^x\|_* \|\bsx_{i,t} -\tilde{\bsx}_{i,t}  \|.
\end{flalign}
Thus,  $(i)$ is derived by some simplification. Similarly, $(ii)$ can be obtained. \hfill$\square$

\begin{lemma} \label{lemma x's diffience}
 Let Assumptions \ref{graph assumption 1}, \ref{assu schostic gradient} and \ref{assu DRX Y noexpansive} hold. Then, we have that for $T\geq 2$,
\begin{flalign}
&(i)\ \mathbb{E} \left[ \sum_{t=1}^T \sum_{i=1}^n \| \bsx_{i,t} - \bsx_{avg,t}\|\right]\leq  \sum_{i=1}^n\| \bsx_{i,1} - \bsx_{avg,1}\| \nonumber \\
& \quad \quad  +  \frac{ n \Gamma}{1-\sigma} \sum_{j=1}^n \|\bsx_{j,1}\| +\frac{ \Gamma }{ 1-\sigma} E_1 \sum_{t=1}^{T-1} \alpha_t, \nonumber \\
&(ii)\ \mathbb{E} \left[  \sum_{t=1}^T \sum_{i=1}^n \| \bsy_{i,t} - \bsy_{avg,t}\|\right]\leq   \sum_{i=1}^n\| \bsy_{i,1} - \bsy_{avg,1}\| \nonumber \\
& \quad \quad   +   \frac{ n \Gamma}{1-\sigma} \sum_{j=1}^n \|\bsy_{j,1}\| +\frac{ \Gamma }{1-\sigma} E_2 \sum_{t=1}^{T-1} \eta_t
\end{flalign}
where $E_1={n^2 L_X}/ { \varrho_x } $ and $E_2= {n^2 L_Y}/ { \varrho_y } $.
\end{lemma}
\noindent{\em Proof:} See Appendix \ref{appendices A}.


The bounds in Lemma \ref{lemma xx yy} describe  the differences between $\tilde{\bsx}_{i,t}$ and $ {\bsx}_{i,t}$ ($\tilde{\bsy}_{i,t}$ and $ {\bsy}_{i,t}$) in Algorithm \ref{algorithm 1}. Lemma \ref{lemma x's diffience} describes the consistency penalty incurred for the  decision disagreement between agents.
Next, we exactly establish the upper-bound of $\textbf{ESP-Regret}_{d}^j (T)$ for Algorithm \ref{algorithm 1}.
\vspace{0.2cm}
\begin{theorem}\label{theorem 1}
Let Assumptions \ref{graph assumption 1}-\ref{assu DRX Y noexpansive}
hold. Then, for $T \geq 2$ and any $j \in \mathcal{V}$, we obtain
\begin{flalign} \label{theorem 1 equation}
&\textbf{ESP-Regret}_d^j(T)   \leq I_1 +\frac{\Gamma}{1-\sigma} I_2 +  \left(I_3 +\frac{\Gamma}{1-\sigma} I_4 \right) \sum\limits_{t=1}^{T} \alpha_t  \nonumber \\
&+   \left(I_5 +\frac{\Gamma}{1-\sigma} I_6 \right) \sum\limits_{t=1}^{T} \eta_t + \frac{n R_X}{\alpha_T} + \frac{n R_Y}{\eta_T} + n K_X \sumT \frac{1}{\alpha_t} \cdot\nonumber \\
&  \| \bsx_{t+1}^*-B_t \bsx_t^*\| +  n K_Y \sumT \frac{1}{\eta_t}  \| \bsy_{t+1}^*-C_t \bsy_t^*\|
\end{flalign}
where
\allowdisplaybreaks{\begin{flalign*}
 &I_1=(n+2)\sum_{i=1}^n \left(L_X\| \bsx_{i,1} - \bsx_{avg,1}\| + L_Y \| \bsy_{i,1} - \bsy_{avg,1}\| \right),\\
  &I_2= n (n+2) \sum_{i=1}^n \left(L_X \|\bsx_{i,1}\|+L_Y \|\bsy_{i,1}\|\right),\\
&I_3 =\frac{n L_X^2}{\varrho_x}, I_4 =\frac{n^2(n+2) L_X^2}{\varrho_x},\\
&I_5 =\frac{n L_Y^2}{\varrho_y}, I_6 =\frac{n^2(n+2) L_Y^2}{\varrho_y},\\
& R_X=\max_{\bsx, \bsu \in \bsX}\Psi_{\mathcal{R}}^x (\bsx, \bsu), R_Y=\max_{\bsy, \boldsymbol{z} \in \bsY} \Psi_{\mathcal{R}}^y (\bsy, \boldsymbol{z}).
\end{flalign*}}
\end{theorem}
\vspace{0.2cm}


\noindent{\em Proof:} Define  the auxiliary regret: \emph{dynamic partial regret}
$\textbf{P-Reg}_d^x (T) = \sum_{t=1}^T \sum_{i=1}^n [f_{i,t} (\boldsymbol{x}_{i,t},\boldsymbol{y}_{i,t}) -  f_{i,t} (\boldsymbol{x}_t^*,\boldsymbol{y}_{i,t})],$  \\
$\textbf{P-Reg}_d^y (T) = \sum_{t=1}^T \sum_{i=1}^n [f_{i,t} (\boldsymbol{x}_{i,t},\boldsymbol{y}_t^*) - f_{i,t}  (\boldsymbol{x}_{i,t},\boldsymbol{y}_{i,t})].$

From (\ref{fun lipschitz}), it can be obtained that
\begin{flalign}\label{proof T1 eq1}
&\sum_{i=1}^n |f_{i,t} (\boldsymbol{x}_{j,t},\boldsymbol{y}_{j,t}) -  f_{i,t} (\boldsymbol{x}_{i,t},\boldsymbol{y}_{i,t}) |\nonumber \\
&\leq \sum_{i=1}^n \left(L_X \| \boldsymbol{x}_{i,t}-\boldsymbol{x}_{j,t}\| +L_Y\|\boldsymbol{y}_{i,t}-\boldsymbol{y}_{j,t}\|\right)\nonumber \\
&\leq (n+1) \sum_{i=1}^n \left(L_X \| \boldsymbol{x}_{i,t}-\boldsymbol{x}_{avg,t}\| + L_Y\|\boldsymbol{y}_{i,t}-\boldsymbol{y}_{avg,t}\|\right).
\end{flalign}

Along (\ref{proof T1 eq1}), we have that
\begin{flalign} \label{proof T1 eq2}
&\textbf{ESP-Regret}_{d}^j (T)  \nonumber \\
& \overset{(a)}{\leq}    \mathbb{E} \left\{  \sumT \sumn \left| f_{i,t} (\bsx_{j,t},\bsy_{j,t}) - f_{i,t} (\bsx_{i,t},\bsy_{i,t}) \right| \right\}  \nonumber \\
&\quad +\left|\mathbb{E} \left\{  \sumT \sumn \left[ f_{i,t} (\bsx_{i,t},\bsy_{i,t})  -  f_{i,t} (\bsx_t^*,\bsy_t^*) \right]  \right\} \right| \nonumber \\
& \overset{(b)}{\leq}    (n+2) \mathbb{E} \left[ \sumT \sumn (L_X \| \boldsymbol{x}_{i,t}-\boldsymbol{x}_{avg,t}\| \right. \nonumber \\
&\quad\left. + L_Y\|\boldsymbol{y}_{i,t}-\boldsymbol{y}_{avg,t}\|) \right] + \Omega_{PR}^x(T) + \Omega_{PR}^y(T)
\end{flalign}
where $\Omega_{PR}^x(T)$ and $\Omega_{PR}^y(T)$ represents two non-negative upper bounds that satisfy $ \mathbb{E} [\textbf{P-Reg}_d^x (T)]\leq \Omega_{PR}^x(T)$ and $ \mathbb{E} [\textbf{P-Reg}_d^y (T)]\leq \Omega_{PR}^y(T)$, respectively, $(a)$ follows triangle inequality and the fact that $|\mathbb{E}[s_1]| \leq\mathbb{E}[ |s_1|], s_1\in\mathbb{R}$, and $(b)$ follows Lemma 2 of \cite{10239326}.

Next, we focus on $\Omega_{PR}^x(T)$ and $ \Omega_{PR}^y(T) $ on the RHS of (\ref{proof T1 eq2}). On the one hand, from the convexity of the loss function over $\bsx$, we obtain that
\begin{flalign}\label{proof T1 eq3}
& \textbf{P-Reg}_{d}^x (T)\nonumber \\
 &  \leq  \sumT \sumn \langle \widetilde{\nabla}_{i,t}^x,  \bsx_{i,t}- \bsx_t^* \rangle    +  \sumT \sumn \langle \boldsymbol{g}_{i,t}^x - \widetilde{\nabla}_{i,t}^x,  \bsx_{i,t}- \bsx_t^* \rangle  \nonumber \\
 &\overset{(c)}{\leq} \sumT \sumn\left\| \widetilde{\nabla}_{i,t}^x \right\|_*  \| \bsx_{i,t}-\tilde{\bsx}_{i,t} \| + \sumT \sumn\langle \widetilde{\nabla}_{i,t}^x, \tilde{\bsx}_{i,t}-\bsx_t^* \rangle  \nonumber \\
 &\quad +  \sumT \sumn \langle  \boldsymbol{g}_{i,t}^x - \widetilde{\nabla}_{i,t}^x,  \bsx_{i,t}- \bsx_t^* \rangle  \nonumber \\
 &\overset{(d)}{\leq} \sumT \sumn \frac{\alpha_t}{\varrho_x}\left\| \widetilde{\nabla}_{i,t}^x \right\|_*^2 + \sumT \sumn\frac{1}{\alpha_t} [ \Psi_{\mathcal{R}}^x( \bsx_t^*, \bsx_{i,t}) \nonumber \\
 &\quad-\Psi_{\mathcal{R}}^x( \bsx_t^*, \tilde{\bsx}_{i,t})]+  \sumT \sumn \langle  \boldsymbol{g}_{i,t}^x - \widetilde{\nabla}_{i,t}^x,  \bsx_{i,t}- \bsx_t^* \rangle
\end{flalign}
where $(c)$ is obtained by using Cauchy-Schwarz inequality and $(d)$ is derived by using the Lemma \ref{lemma xx yy}, (\ref{lem xx yy eq02}), and the fact that $\Psi_{\mathcal{R}}^x(\tilde{\bsx}_{i,t}, \bsx_{i,t})\geq 0$.

Now, we turn our attention to the term $\Psi_{\mathcal{R}}^x( \bsx_t^*, \bsx_{i,t}) -\Psi_{\mathcal{R}}^x( \bsx_t^*, \tilde{\bsx}_{i,t})$ in (\ref{proof T1 eq3}). By adding and subtracting the terms $\Psi_{\mathcal{R}}^x( \bsx_{t+1}^*, \bsx_{i,t+1})$ and $\Psi_{\mathcal{R}}^x(B_t \bsx_t^*, \bsx_{i,t+1})$, we yield that
\begin{flalign}\label{proof T1 eq4}
 & \Psi_{\mathcal{R}}^x( \bsx_t^*, \bsx_{i,t}) -\Psi_{\mathcal{R}}^x( \bsx_t^*, \tilde{\bsx}_{i,t})\nonumber \\
 &=  \Psi_{\mathcal{R}}^x( \bsx_t^*, \bsx_{i,t}) -\Psi_{\mathcal{R}}^x( \bsx_{t+1}^*, \bsx_{i,t+1})\nonumber \\
 &\quad  +\Psi_{\mathcal{R}}^x( \bsx_{t+1}^*, \bsx_{i,t+1}) -  \Psi_{\mathcal{R}}^x(B_t \bsx_t^*, \bsx_{i,t+1}) \nonumber \\
 &\quad +  \Psi_{\mathcal{R}}^x(B_t  \bsx_t^*, \bsx_{i,t+1}) -\Psi_{\mathcal{R}}^x( \bsx_t^*, \tilde{\bsx}_{i,t}).
\end{flalign}
Further, these terms on the RHS of (\ref{proof T1 eq4}) follow that
\begin{flalign}\label{proof T1 eq5}
  1)\ & \sum_{t=1}^T \sum_{i=1}^n \frac{1}{\alpha_t} \left[  \Psi_{\mathcal{R}}^x (\bsx_t^*, \bsx_{i,t}) -\Psi_{\mathcal{R}}^x( \bsx_{t+1}^*, \bsx_{i,t+1}) \right] \nonumber \\
   &\leq  \frac{n R_X}{\alpha_1}  +  n R_X \sum_{t=2}^T  \left(\frac{1}{\alpha_t} -\frac{1}{\alpha_{t-1}} \right)\nonumber \\
   & \leq \frac{n R_X}{\alpha_T},  \\
  2)\ & \Psi_{\mathcal{R}}^x( \bsx_{t+1}^*, \bsx_{i,t+1}) -  \Psi_{\mathcal{R}}^x(B_t  \bsx_t^*, \bsx_{i,t+1})\nonumber \\
  &\leq K_X    \|  \bsx_{t+1}^*- B_t \bsx_t^*\|,  \\
  %
  3)\ & \sum_{i=1}^n \frac{1}{\alpha_t}\left[  \Psi_{\mathcal{R}}^x (B_t  \bsx_t^*, \bsx_{i,t+1}) -\Psi_{\mathcal{R}}^x( \bsx_t^*, \tilde{\bsx}_{i,t})\right]\nonumber \\
  & \leq \sum_{i=1}^n \frac{1}{\alpha_t} \sum_{j=1}^n [A_t]_{ij}  \left[  \Psi_{\mathcal{R}}^x (B_t  \bsx_t^*, \bss_{j,t}^x) -\Psi_{\mathcal{R}}^x( \bsx_t^*, \tilde{\bsx}_{i,t})\right]\nonumber \\
  & \leq 0
\end{flalign}
where $1)$ follows the facts $\frac{1}{\alpha_t}-\frac{1}{\alpha_{t-1}}\geq 0$ and $\Psi_{\mathcal{R}}^x(\bsx_1, \bsx_2)\geq 0, \forall \bsx_1, \bsx_2 \in \bsX$,  $2)$ follows  Assumption \ref{assu DR continous}, and $3)$ is obtained by using the double stochasticity of $A_t$ and Assumption \ref{assu DRX Y noexpansive}.

Substituting them into $\Psi_{\mathcal{R}}^x( \bsx_t^*, \bsx_{i,t}) -\Psi_{\mathcal{R}}^x( \bsx_t^*, \tilde{\bsx}_{i,t})$ and summing it over $i\in[n]$ and $t\in[T]$, it follows that
\begin{flalign}\label{proof T1 eq6}
&\sumT \sumn \frac{1}{\alpha_t} \left[ \Psi_{\mathcal{R}}^x( \bsx_t^*, \bsx_{i,t}) -\Psi_{\mathcal{R}}^x( \bsx_t^*, \tilde{\bsx}_{i,t}) \right] \nonumber \\
 & \leq  \frac{n R_X}{\alpha_T} + n K_X \sumT \frac{1}{\alpha_t}  \| \bsx_{t+1}^*-B_t \bsx_t^*\|.
\end{flalign}

Taking expectation operation from the term $ \langle \boldsymbol{g}_{i,t}^x - \widetilde{\nabla}_{i,t}^x,  \bsx_{i,t}- \bsx_t^* \rangle$ in (\ref{proof T1 eq3}), we yield that
$
\mathbb{E} [ \langle \boldsymbol{g}_{i,t}^x - \widetilde{\nabla}_{i,t}^x,  \bsx_{i,t}- \bsx_t^* \rangle ]
  = \mathbb{E} \{ \mathbb{E} [ \langle \boldsymbol{g}_{i,t}^x - \widetilde{\nabla}_{i,t}^x,  \bsx_{i,t}- \bsx_t^* \rangle \big| \mathcal{F}_{t-1} ] \}
 = \mathbb{E} \{ \langle \mathbb{E} [ \boldsymbol{g}_{i,t}^x - \widetilde{\nabla}_{i,t}^x \big| \mathcal{F}_{t-1} ],  \bsx_{i,t}- \bsx_t^* \rangle \}
 = 0.
$
Further, based on (\ref{proof T1 eq6}), and the fact that $\mathbb{E} [ \|\widetilde{\nabla}_{i,t}^{x} \|_*^2 ]=  \mathbb{E} \{  \mathbb{E} [ \|\widetilde{\nabla}_{i,t}^{x} \|_*^2 |  \mathcal{F}_{t-1} ] \} \leq L_X^2$, we get
\begin{flalign}\label{proof T1 eq7}
 \mathbb{E} [\textbf{P-Reg}_{d}^x (T)] & \leq  \frac{n L_X^2}{\varrho_x}\sumT \alpha_t+ n K_X \sumT \frac{1}{\alpha_t}  \| \bsx_{t+1}^*-B_t \bsx_t^*\|\nonumber \\
&\quad + n R_X/\alpha_T.
\end{flalign}

On the other hand,  it can be similarly obtained that
$
 \mathbb{E}[\textbf{P-Reg}_{d}^y (T)] \leq  \frac{n L_Y^2}{\varrho_y}\sum_{t=1}^T \eta_t + n K_Y \sum_{t=1}^T \frac{1}{\eta_t}  \| \bsy_{t+1}^*-C_t \bsy_t^*\|+ n R_Y/\eta_T.
$
Finally, by combining the inequality, (\ref{proof T1 eq2}), (\ref{proof T1 eq7}), and Lemma \ref{lemma x's diffience}, Theorem \ref{theorem 1} is established.
\hfill$\square$

In Theorem \ref{theorem 1}, the convergence result indicates that the dynamic regret bound is significantly influenced by the tunning step sizes $\alpha_t, \eta_t$. Therefore, we further explore the exact convergence rate of Algorithm \ref{algorithm 1} in Corollary \ref{corollary 1}.

\begin{corollary} \label{corollary 1}
 Suppose that the conditions required in Theorem \ref{theorem 1} and $V_T = o(T)$ hold. Taking
$
 \alpha_t=\frac{1}{\epsilon_1} t^{- \gamma_1},  \eta_t=\frac{1}{\epsilon_2 } t^{- \gamma_2}, \epsilon_1, \epsilon_2 > 0,  \gamma_1,  \gamma_2 \in (0, 1),
$
 we yield for $T \geq 2$ and $j \in \mathcal{V}$ that,
\begin{flalign}\label{corollary 1 equation}
&\textbf{ESP-Regret}_d^j(T) \nonumber\\
&\leq \mathcal{O}\left(\max \left\{\left(1+\frac{\Gamma}{1-\sigma}\right) T^{\theta_1},  T^{\theta_2} (1+V_T ) \right\}\right)
\end{flalign}
where $\theta_1=\max \{1-\gamma_1, 1-\gamma_2 \}, \theta_2=\max \{\gamma_1, \gamma_2 \}$.
\end{corollary}

\noindent{\em Proof:} By substituting this step size set in Corollary \ref{corollary 1} into (\ref{theorem 1 equation}), we yield that
\begin{flalign}
&\textbf{ESP-Regret}_d^j(T)\nonumber\\
 &\leq I_1 + \frac{\Gamma I_{2}}{1-\sigma} + \epsilon_1^{-1} \left(I_3+ \frac{\Gamma I_4}{ 1-\sigma}\right) \sum_{t=1}^{T} t^{-\gamma_1} \nonumber\\
  & \quad +  \epsilon_2^{-1} \left(I_5+\frac{\Gamma I_6}{1-\sigma}\right) \sum_{t=1}^{T} t^{-\gamma_2} + n R_X\epsilon_1 T^{\gamma_1}+ n R_Y \epsilon_2 T^{\gamma_2} \nonumber\\
 &\quad + n K_X \epsilon_1 T^{\gamma_1} V_T^x
 +  n K_Y  \epsilon_2 T^{\gamma_2} V_T^y\nonumber\\
&\leq \mathcal{O} \left(\left(1+ \frac{\Gamma}{1-\sigma}\right) (T^{1-\gamma_1}+T^{1-\gamma_2}) +T^{\gamma_1}(1+V_T^x) \right. \nonumber\\
 &\quad \left. +T^{\gamma_2}(1+V_T^y) \right)
\end{flalign}
where the last inequality follows that
\begin{flalign*}
\sum_{t=1}^T \frac{1}{t^{\gamma_1}} =1+\sum_{t=2}^T \frac{1}{t^{\gamma_1}} \leq 1+ \int_1^T \frac{1}{t^{\gamma_1}} dt \leq \frac{1}{1-\gamma_1} T^{1-\gamma_1}.
\end{flalign*}
Based on the definitions of $\theta_1, \theta_2$, and $V_T$, (\ref{corollary 1 equation}) is gotten.
\hfill$\square$

\begin{remark}
From Corollary \ref{corollary 1},  it can be known that the different choices of coefficients $\gamma_1$ and $\gamma_2$ correspond to different orders of $\textbf{ESP-Regret}_d^j(T)$ over $T$. Specially, when $\gamma_{1}= \gamma_{2}=1/2$ holds,  the regret bound $ \mathcal{O}(  \sqrt{T} (1+V_T ) )$  is obtained; when $\gamma_{1}= \gamma_{2}= 1/2-\log_T {\sqrt{1+V_T}}$ holds, the regret bound $\mathcal{O}( \sqrt{ T (1+V_T)} )$  is obtained. The latter bound implemented by our algorithm matches the optimal convergence performance of the centralized counterpart, provided that the knowledge of $V_T$ is known \cite{zhang2018adaptive}.
\end{remark}

\begin{remark}
Compared with the distributed online subgradient saddle point algorithm proposed in \cite{10239326}, the dynamic regret bound (\ref{corollary 1 equation}) achieved by Algorithm \ref{algorithm 1} are more general due to its non-Euclidean distance metrics and the path-variation $V_T$ with predictive mappings. Further, benefitting from the path-variations defined in (\ref{assumption Path Variation equation}), the regret bound enables potential performance gains through finding appropriate mapping $B_t$ and $C_t$ to achieve a low order of $V_T$ with respect to $T$.
\end{remark}

\section{Multiple Consensus Iterations} \label{s4}
In the section, we investigate a multiple consensus version of Algorithm DOSMD-CCO (Multi-DOSMD-CCO) and aim to improve the consensus process among agents by employing a multiple consensus technique. By conducting multiple consensus steps in each iteration, each agent can effectively integrate local information from agents further away to accelerate global consensus and improve optimization efficiency\cite{eshraghi2020distributed,jakovetic2014fast}.  In addition, for limited communication networks, such as network delays and noise, multiple consensus has better robustness than one consensus.
Taking the collaborative search and rescue of robot groups as an example, multiple consensus technique allows each robot to receive more information from its companions in each iteration, such as their position changes and task status, which enables them to better adjust search strategies and improve completion efficiency in a real-time environment.



 Different from the step 6 of Algorithm \ref{algorithm 1}, agent $i$ carries out $K_t \in \mathbb{Z}_+$ times of consensus operations at each time $t$. This implies that each agent can receive the information of other agents that are $K_t$ hops away. Let $\bss_{i,t}^{x,1}= \bss_{i,t}^x$ and $\bss_{i,t}^{y,1}= \bss_{i,t}^y, i\in[n]$. Then, agent $i$ receives $\boldsymbol{s}_{j,t}^{x,1}$ and $\boldsymbol{s}_{j,t}^{y,1}$ from $j \in \mathcal{N}_{i}^{\text {in }}(t)$, and updates over $k=1$ to $K_t$ that
\begin{flalign} \label{Multiple Consensus Iterations}
\bss_{i,t}^{x,k+1}&= \sum_{j \in \mathcal{N}_{i}^{in}(t)}{[A_t]_{ij} \bss_{j,t}^{x,k}}, \nonumber \\
\bss_{i,t}^{y,k+1}&= \sum_{j \in \mathcal{N}_{i}^{in}(t)}[A_t]_{ij} \bss_{j,t}^{y,k}.
\end{flalign}
Then, let  $\bsx_{i,t+1}=\bss_{i,t}^{x,K_t+1},\ \bsy_{i,t+1}=\bss_{i,t}^{y,K_t+1}.$ Based on this consensus operations, the Multi-DOSMD-CCO algorithm  is presented in Algorithm \ref{algorithm 2}.

\begin{algorithm}[!t]
	\renewcommand{\algorithmicrequire}{\textbf{Initialize:} }
	\caption{  Multi-DOSMD-CCO algorithm.}
	\label{algorithm 2}
	\begin{algorithmic}[1]
		\REQUIRE Initial  decisions $\boldsymbol{x}_{i,1} \in \boldsymbol{X} ,\boldsymbol{y}_{i,1} \in \boldsymbol{Y},$ the parameters $\alpha_t,\eta_t>0, K_t\in \mathbb{Z}_+$, and the mappings $B_t, C_t$. Set $\bss_{j,t}^{x,1}= \bss_{j,t}^x$ and $ \bss_{j,t}^{y,1}= \bss_{j,t}^y$.
\setlength{\parskip}{0.4em}
\FOR {$t=1,2,\cdots,T$}
             \FOR {$i \in \mathcal{V}$ in parallel}

             \STATE Agent $i$ gets the stochastic gradients $\widetilde{\nabla}_{i,t}^x, \widetilde{\nabla}_{i,t}^y$, and computes, respectively,
 \setlength{\parskip}{0.2em}

\begin{center}
$\nabla \mathcal{R}^x (\bsz_{i,t})= \nabla \mathcal{R}^x (\bsx_{i,t}) -\alpha_t \widetilde{\nabla}_{i,t}^x,$

$\nabla \mathcal{R}^y (\bsv_{i,t})= \nabla \mathcal{R}^y (\bsy_{i,t}) +\eta_t \widetilde{\nabla}_{i,t}^y.$
\end{center}

\STATE Executes the Bregman projections, respectively,
\begin{center}
$\tilde{\bsx}_{i,t}= \underset{\bsx\in\bsX}{\arg\min} \ \Psi_{\mathcal{R}}^x (\bsx, \bsz_{i,t}),$

$\tilde{\bsy}_{i,t}=\underset{\bsy\in\bsY}{\arg\min}\ \Psi_{\mathcal{R}}^y (\bsy, \bsv_{i,t}).$
\end{center}

\STATE Runs the decision corrections using prediction mappings, i.e.,

\begin{center}
 $\bss_{i,t}^x= B_t \tilde{\bsx}_{i,t},\ \bss_{i,t}^y= C_t \tilde{\bsy}_{i,t}.$		
\end{center}
		
\STATE Runs $K_t$ consensus steps:
        \FOR {$k =1,2,\cdots,K_t$}

            \STATE Agent $i$ receives $\boldsymbol{s}_{j,t}^{x,k}, \boldsymbol{s}_{j,t}^{y,k}$ from its in-neighbors, and updates

              \begin{center}
              $\bss_{i,t}^{x,k+1}= \sum\limits_{j \in \mathcal{N}_{i}^{in}(t)}{[A_t]_{ij} \bss_{j,t}^{x,k}}$,

           		 $\bss_{i,t}^{y,k+1}= \sum\limits_{j \in \mathcal{N}_{i}^{in}(t)}[A_t]_{ij} \bss_{j,t}^{y,k}$.
              \end{center}

		\ENDFOR

         \STATE Let  $\bsx_{i,t+1}=\bss_{i,t}^{x,K_t+1}$ and $ \bsy_{i,t+1}=\bss_{i,t}^{y,K_t+1}.$
		\ENDFOR
        \ENDFOR
	\end{algorithmic}
\end{algorithm}

 By observing Algorithm \ref{algorithm 2}, it is not hard to note that in a theoretical sense, (\ref{Multiple Consensus Iterations}) are equivalent to
\begin{flalign}
\bsx_{i,t+1}&=  \sum_{j =1}^n {[A_t]_{ij} \bss_{j,t}^{x,K_t}}= \sum_{j=1}^n {\left[(A_t)^{K_t}\right]_{ij} \bss_{j,t}^x}, \nonumber \\
\boldsymbol{y}_{i,t+1}&=  \sum_{j =1}^n {[A_t]_{ij} \bss_{j,t}^{y,K_t}}= \sum_{j=1}^n {\left[(A_t)^{K_t}\right]_{ij} \bss_{j,t}^y}, \label{MConsensus equ}
\end{flalign}

Now, we turn our attention to the convergence performance of Algorithm Multi-DOSMD-CCO. To ensure the feasibility of the multiple consensus method at any time $t$, the following assumption is made, which is set similarly in \cite{xiong2022event, yuan2020distributed}.
\begin{assumption}\label{graph assumption 2}
$\mathcal{G}_t$ is strongly connected for time $t \in [T]$.
\end{assumption}

Next, for $t\geq s \geq 1$, we introduce the transition matrix
\begin{flalign} \label{transition matrix 2}
\Phi^K(t,s)\triangleq (A_t)^{K_t}(A_{t-1})^{K_{t-1}}\cdots (A_{s})^{K_{s}}
\end{flalign}
and set $ \Phi^K(t,t+1)=I$.
Based on (\ref{transition matrix 2}), Assumption \ref{graph assumption 2} and the convergence property of $\Phi(t,s)$ stated in Lemma \ref{lemma-graph},
we derive the following condition:
 \begin{flalign} \label{phi-K-b2}
\left|[\Phi^K(t, s)]_{i j}-\frac{1}{n}\right| \leq \Gamma_1 \sigma_1^{\sum_{p=s}^tK_p-1}
\end{flalign}
where $\Gamma_1=(1-\zeta / 4 n^{2})^{-1}$ and $\sigma_1=(1-\zeta / 4 n^{2})$.

Now, we establish the bound of the consensus errors.

\begin{lemma} \label{lemma x's diffience2}
Suppose that Assumptions \ref{assu DRX Y noexpansive} and \ref{graph assumption 2} hold. Then, we have that for $T\geq 2$,
\begin{flalign}
&(i)\ \mathbb{E} \left[ \sum_{t=1}^T \sum_{i=1}^n \| \bsx_{i,t} - \bsx_{avg,t}\|\right]\leq   \sum_{i=1}^n\| \bsx_{i,1} - \bsx_{avg,1}\|  \nonumber \\
&\quad \quad +\frac{ n \Gamma_1 \sigma_1^{\underline{K}-1}}{1-\sigma_1^{\underline{K}}} \sum_{j=1}^n \|\bsx_{j,1}\|+ \frac{ \Gamma_1 \sigma_1^{\underline{K}-1}}{ 1-\sigma_1^{\underline{K}}} E_1 \sum_{t=1}^{T-1} \alpha_t,  \\
& (ii)\  \mathbb{E} \left[  \sum_{t=1}^T \sum_{i=1}^n\| \bsy_{i,t} - \bsy_{avg,t}\|\right]  \leq \sum_{i=1}^n\| \bsy_{i,1} - \bsy_{avg,1}\|\nonumber \\
& \quad \quad +\frac{ n \Gamma_1 \sigma_1^{\underline{K}-1}}{1-\sigma_1^{\underline{K}}} \sum_{j=1}^n \|\bsy_{j,1}\|
+ \frac{\Gamma_1  \sigma_1^{\underline{K}-1}}{  1-\sigma_1^{\underline{K}} } E_2 \sum_{t=1}^{T-1} \eta_t
\end{flalign}
where $\underline{K}=\min_{t\in[T]} \{ K_t\}.$
\end{lemma}
\noindent{\em Proof:} See Appendix \ref{appendices A}.

\begin{theorem}\label{theorem 2}
Let Assumptions \ref{assu schostic gradient}-\ref{graph assumption 2}
hold and $\{\boldsymbol{x}_{i,t},\boldsymbol{y}_{i,t}\}, i\in[n]$ be the decision sequence obtained from Algorithm Multi-DOSMD-CCO. Then,  we obtain for $T \geq 2$ and $j \in \mathcal{V}$ that
\begin{flalign}\label{theorem 2 equation}
&\textbf{ESP-Regret}_d^j(T) \nonumber \\
 &\leq I_1 + \frac{{\Gamma_1 \sigma_1}^{\underline{K}-1}}{1-{\sigma_1}^{\underline{K}}}I_2 + \left(I_{3}+\frac{{\Gamma_1 \sigma_1}^{\underline{K}-1}}{1-{\sigma_1}^{\underline{K}}}I_4\right) \sum\limits_{t=1}^{T} \alpha_t \nonumber \\
 &+   \left(I_{5}+\frac{\Gamma_1 {\sigma_1}^{\underline{K}-1}}{1-{\sigma_1}^{\underline{K}}}I_6\right) \sum\limits_{t=1}^{T} \eta_t + \frac{n R_X}{\alpha_T} + \frac{n R_Y}{\eta_T} +\Xi_T
\end{flalign}
where $\Xi_T= n K_X \sum_{t=1}^T \frac{1}{\alpha_t}  \| \bsx_{t+1}^*-B_t \bsx_t^*\| +  n K_Y \sum_{t=1}^T\frac{1}{\eta_t} $ $\| \bsy_{t+1}^*-C_t \bsy_t^*\|.$
\end{theorem}
\vspace{0.15cm}

\noindent{\em Proof:} Based on the inequality
\begin{flalign}
    &\sum_{i=1}^n \frac{1}{\alpha_t} \left[  \Psi_{\mathcal{R}}^x(B_t  \bsx_t^*, \bsx_{i,t+1}) -\Psi_{\mathcal{R}}^x( \bsx_t^*, \tilde{\bsx}_{i,t})\right]\nonumber \\
    &\leq \sum_{i=1}^n \frac{1}{\alpha_t} \sum_{j=1}^n [(A_t)^{K_t} ]_{ij}  \left[ \Psi_{\mathcal{R}}^x (B_t  \bsx_t^*, \bss_{j,t}^x) -\Psi_{\mathcal{R}}^x( \bsx_t^*, \tilde{\bsx}_{i,t})\right]\nonumber \\
    &\leq 0,
\end{flalign}
the bounds of  $\mathbb{E} [\textbf{P-Reg}_{d}^{x[y]} (T)]$ in (\ref{proof T1 eq7}) still hold for Algorithm Multi-DOSMD-CCO. Finally, by combining the bounds, (\ref{proof T1 eq2}), and Lemma \ref{lemma x's diffience2}, (\ref{theorem 2 equation}) is derived. \hfill$\square$

Theorem \ref{theorem 2} strictly characterises the upper-bound of $\textbf{ESP-Regret}_d^j(T)$ for Algorithm Multi-DOSMD-CCO and shows the effects of the multiple consensus parameter $\underline{K}$, the adjustable step sizes $\alpha_t$ and $\eta_t$ on this bound. Further, to clearly reveal the performance gains given from the multiple consensus technique, we explore its exact convergence rate w.r.t. $T$ in Corollary \ref{corollary 2} under the setting of $\alpha_t, \eta_t$.

\begin{corollary} \label{corollary 2}
 Suppose that the conditions required in Theorem \ref{theorem 2}  and $V_T = o(T)$ hold. Taking the same  $\alpha_t, \eta_t$ with Corollary \ref{corollary 1},
 we yield for $T \geq 2$ and $j \in \mathcal{V}$ that,
\begin{flalign}\label{corollary 2 equation}
&\quad \textbf{ESP-Regret}_d^j(T) \nonumber\\
& \leq \mathcal{O}\left( \max\left\{ \left(1+\frac{\Gamma_1{\sigma_1}^{\underline{K}-1}}{1-{\sigma_1}^{\underline{K}}}\right)T^{\theta_1},  T^{\theta_2} (1+V_T ) \right\} \right).
\end{flalign}
\end{corollary}
\vspace{0.15cm}

\noindent{\em Proof:} Similar to the proof of (\ref{corollary 1 equation}), (\ref{corollary 2 equation}) can be established.  \hfill$\square$

\begin{remark} \label{remark coro 2}
The bound in (\ref{corollary 2 equation}) clearly quantifies and reveals the potential performance gain brought from the multiple consensus technique via the parameter $\underline{K}=\min_{t\in[T]} \{ K_t\}$. Based on the fact that $\Gamma=\Gamma_1, \sigma=\sigma_1$ under Assumption \ref{graph assumption 2}, it can be known by comparing the results in Corollaries \ref{corollary 1} and \ref{corollary 2} that Algorithm Multi-DOSMD-CCO can achieve better convergence performance according to the observation that $\frac{\Gamma_1{\sigma_1}^{\underline{K}-1}}{1-{\sigma_1}^{\underline{K}}} \leq \frac{\Gamma_1}{1-\sigma_1}=\frac{\Gamma}{1-\sigma}, \underline{K}\geq1$. And when the larger setting of $\underline{K}$ is satisfied, the performance gain becomes greater, especially for the case when $\sigma_1$ is close to $1$. Naturally, the bound (\ref{corollary 2 equation}) degenerates to the one in Corollary \ref{corollary 1} under Assumption \ref{graph assumption 2} when $K_t=1$ holds.
\end{remark}
\begin{remark}
The performance benefit brought from the multiple consensus technique is often accompanied by the increase of communication overhead, as agents must exchange information more frequently to maintain decision consensus. The challenge lies in balancing the need for multiple consensus rounds to ensure convergence with the communication burden they impose. Excessive communication may reduce the overall optimization efficiency and cause waste of communication resources, especially in large-scale multiagent network. Thus, careful consideration of the trade-off between iteration frequency and communication cost is essential to achieve optimal performance without overwhelming the network.
%
%
\end{remark}

\section{Numerical Simulations} \label{s5}
In this section, we employ the distributed robust tracking problem modeled in (\ref{Simulation example})  as a simulation example to verify the efficacy of the proposed algorithms. By solving the following optimization problem, each agent generates an optimal action in the worst-case sense \cite{ben2009robust, zhang2023communication} to cooperatively track an autonomous moving target.
\begin{flalign} \label{Simulation example}
  \min\limits_{\boldsymbol{x}_t} \max\limits_{\bsy_t}\, &\ \sum_{t=1}^T \sum_{i=1}^n  h_{i,t}(\bsx_t, \bsy_t) + r_{i,t}(\bsx_t, \bsy_t) \nonumber \\
 s.t.& \ \boldsymbol{x}_t \in \boldsymbol{X}, \bsy_t \in \boldsymbol{Y}
\end{flalign}
where $  h_{i,t}(\bsx_t, \bsy_t) = c_{i,1} \langle  \boldsymbol{\pi}_{i,t}, \bsx_t \rangle  +c_{i,2} \| \bsx_t -(\hat{\boldsymbol{\varphi}}_{i,t}+ \boldsymbol{\bsy}_t)  \|_2^2, $ $r_{i,t}(\bsx_t, \bsy_t)= \lambda_{i,1} \|\bsx_t \|_1 -\lambda_{i,2} \|\boldsymbol{\bsy}_t \|_2^2$,
$\boldsymbol{X}:=\{ \boldsymbol{x} | \   [\boldsymbol{x}]_i \in [-6,6], i\in [d] \} $,
$\boldsymbol{Y}:=\{\boldsymbol{\bsy} | \   [\boldsymbol{\bsy}]_i \in [-0.3,0.3], i\in [d] \} $.
For agent $i\in [n]$, $\boldsymbol{\pi}_{i,t}$ and $\hat{\boldsymbol{\varphi}}_{i,t} $ represent the tracking cost vector and the observed value $\hat{\boldsymbol{\varphi}}_{i,t}$ of the target action, respectively. $c_{i,1}$ and $ c_{i,2}$ denote two non-negative coefficients weighing the cost of movement and the effectiveness of tracking. $\lambda_{i,1}$ and $ \lambda_{i,2}$ represent the non-negative regularization parameters.
\begin{figure}[!t]
  \centering
\begin{minipage}{.48\textwidth}
  \centering
\includegraphics[width=8.5cm]{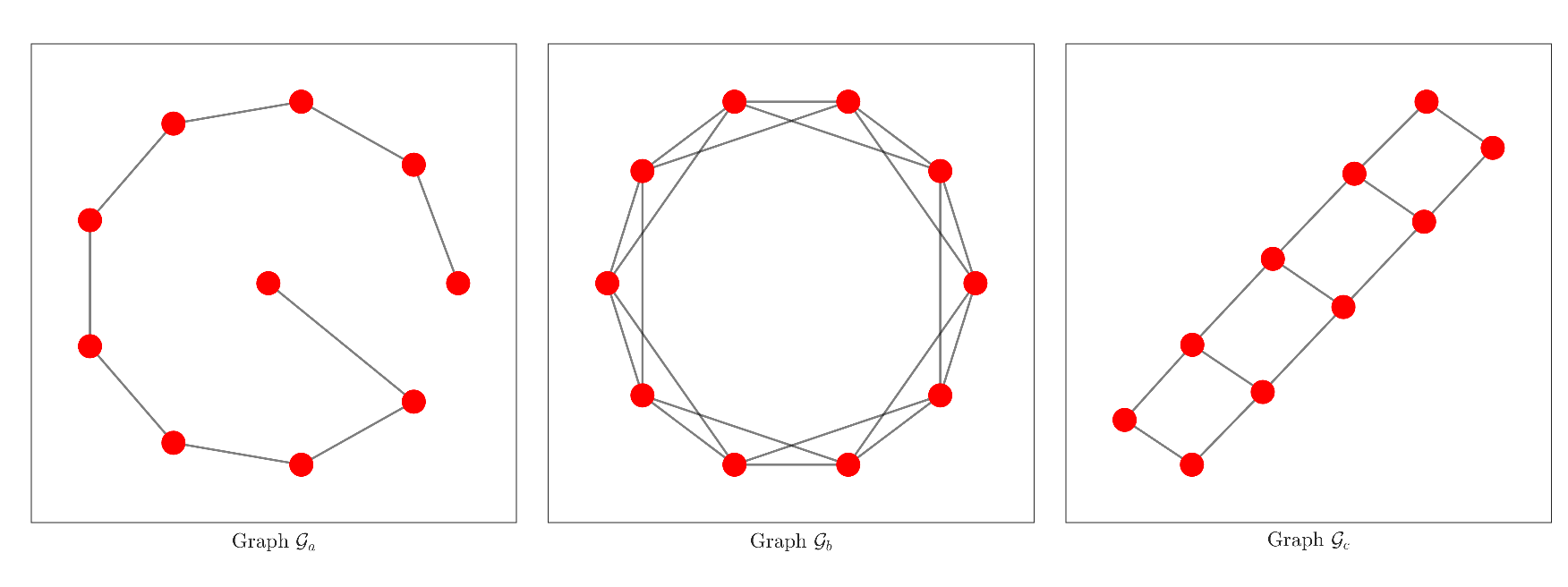}
\caption{The graphs $\mathcal{G}_a, \mathcal{G}_b,$ and $\mathcal{G}_c$.}
\label{f0_three_graph}
\end{minipage}
\begin{minipage}{.48\textwidth}
  \centering
\includegraphics[width=7.2cm]{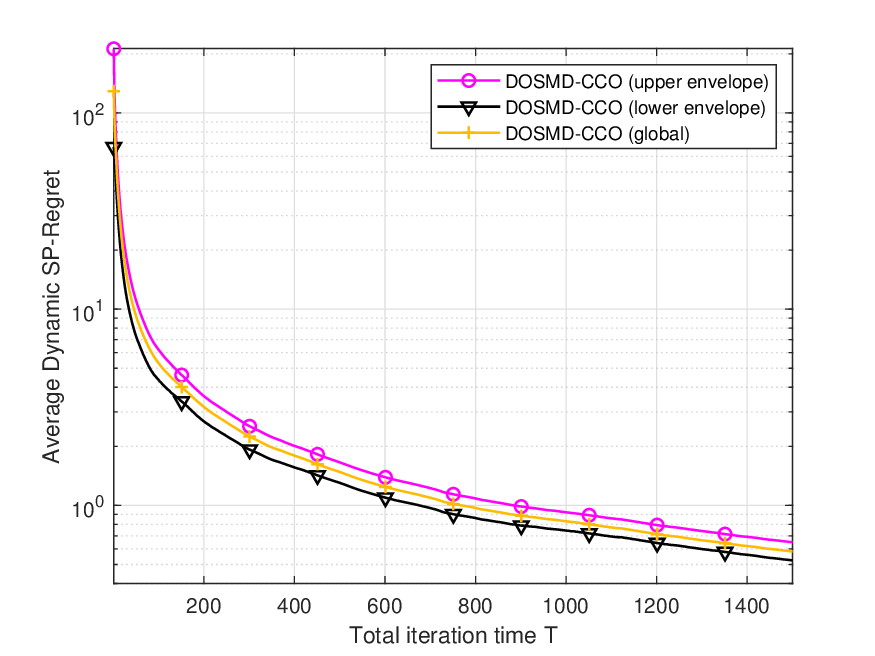}
\caption{The convergence results for Algorithm DOSMD-CCO.}
\label{f1_maxminAVreg}
\end{minipage}
\begin{minipage}{.48\textwidth}
  \centering
\includegraphics[width=7.2cm]{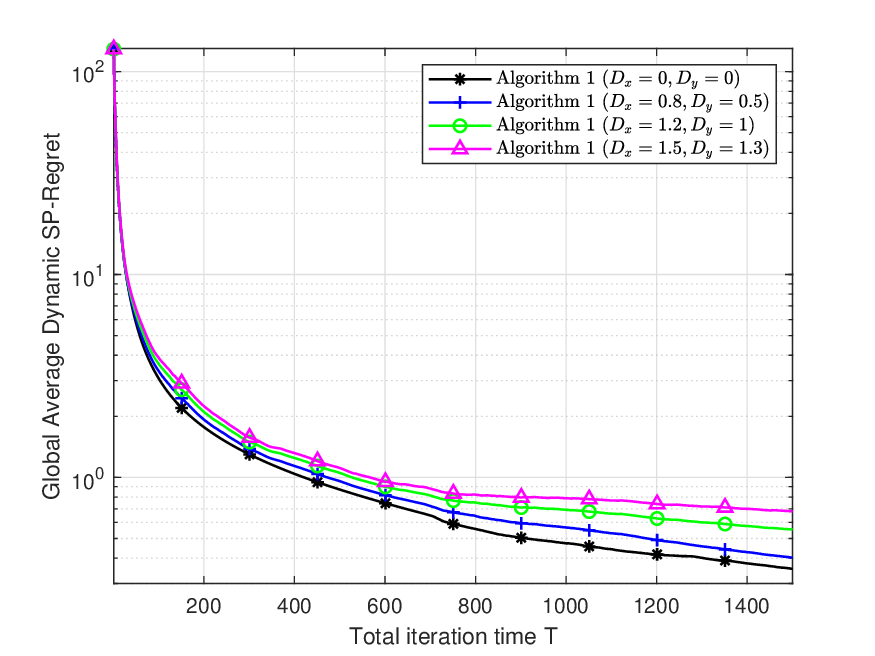}
\caption{The results under different noise magnitudes.}
\label{f3_noise}
\end{minipage}
\vspace{-0.3cm}
\end{figure}

In the following simulations, we set $[\boldsymbol{\pi}_{i,t}]_l=\text{sgn}([\bsx_{i,t}]_l)\cdot [(1-\frac{a_1}{\sqrt{t}}) \boldsymbol{\pi}_i^1 + \frac{a_1}{\sqrt{t}} \boldsymbol{\pi}_{i,t}^2]_l, l\in [d], a_1 \in [0,1)$, where $\boldsymbol{\pi}_i^1$ and $\boldsymbol{\pi}_{i,t}^2$ are two random vectors, and $\text{sgn}([\bsx_{i,t}]_l)$ is a correction term to guarantee that the tracking costs $[\pi_{i,t}]_l \cdot [\bsx_{i,t}]_l$ are always non-negative.
Given the matrix $P_t$ satisfying $\| P_t\|_2 \leq 1$, the target action and its observed variant are set as $\boldsymbol{\varphi}_{t+1}= P_t \boldsymbol{\varphi}_{t} $ and $\hat{\boldsymbol{\varphi}}_{i,t}= \boldsymbol{\varphi}_{t} +\frac{1}{\sqrt{t}} \boldsymbol{\varepsilon}_{i,t}$, respectively, where $\boldsymbol{\varphi}_{1}$ is a random initial vector from $\bsX$ and $\boldsymbol{\varepsilon}_{i,t}$ is a randomly generated observation perturbation.
The gradient noises are generated randomly and independently from the normal distributions $\mathcal{N}(0,D_x 1_d)$ and $\mathcal{N}(0,D_y 1_d)$, respectively.
 let $n=10$, $d=6$, $\mathcal{R}^x(\bsx)=\mathcal{R}^y(\bsx)$ $=\frac{1}{2} \|\bsx \|_2^2, c_{i,1}=0.1, c_{i,2}=0.3, \lambda_{i,1}=0.1$, and $\lambda_{i,2}=0.5$.
Let Algorithms \ref{algorithm 1} and Multi-DOSMD-CCO perform in a time-varying network consisting of three connected graphs $\mathcal{G}_a, \mathcal{G}_b,$ and $\mathcal{G}_c$ shown in Fig. \ref{f0_three_graph}. The network topology switches based on the remainder of round $t$ divided by 3: for 0, $\mathcal{G}_t=$ $\mathcal{G}_a$; for 1, $\mathcal{G}_t=\mathcal{G}_b$; for 2, $\mathcal{G}_t=\mathcal{G}_c$.
 Taking $ \mathcal{M}_j(T):=\textbf{SP-Regret}_d^j (T)/T $, called the average  dynamic  saddle point  regret (ADSPR), as the base metric, we use its upper envelope $\sup_j  \mathcal{M}_j(T)$, its lower envelope $\inf_j  \mathcal{M}_j(T)$, and its global ADSPR $\frac{1}{n}\sum_{j=1}^n  \mathcal{M}_j(T)$ to measure the performance of the proposed algorithm.

Firstly, we study the convergence of Algorithm \ref{algorithm 1} under noise amplitude $D_x=D_y=1$. With  $B_t=C_t=I_d$, the convergence plots around the metric ADSPR, included its upper, lower envelope and global variant, are shown in Fig \ref{f1_maxminAVreg}. From this figure, it is known that Algorithm \ref{algorithm 1} is convergent, which corresponds to Corollary \ref{corollary 1}. Subsequently, by setting different noise magnitudes, we investigate the influence of gradient noises on convergence performance. Fig. \ref{f3_noise} reveals that the optimality of Algorithm \ref{algorithm 1} with exact gradient information is better than these with gradient noises. Moreover, the optimality of the developed algorithm will become worse as the noise magnitudes $D_x$ and $D_y$ is set larger.

\begin{figure}[!t]
  \centering
\begin{minipage}{.48\textwidth}
  \centering
\includegraphics[width=7.5cm]{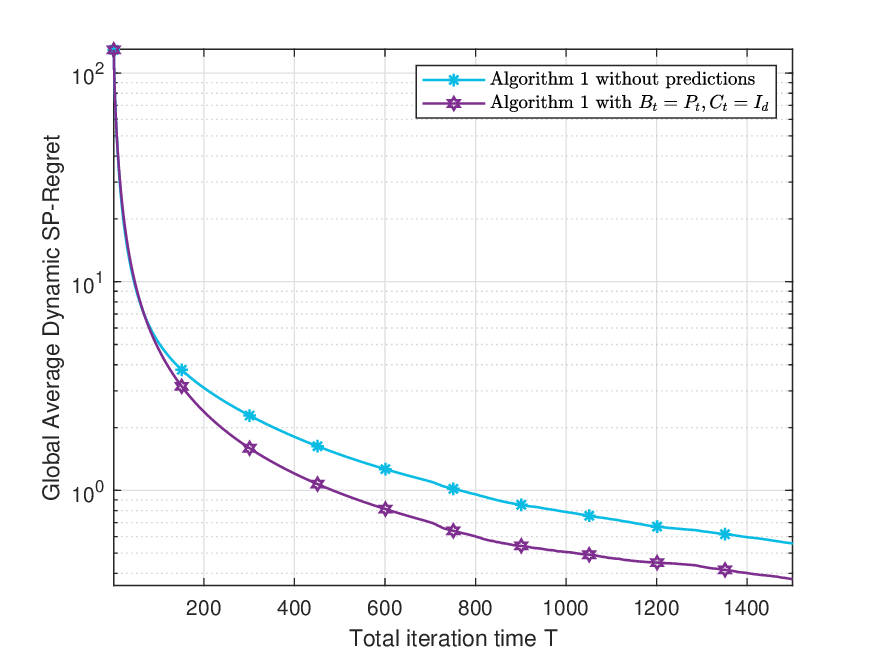}
\caption{The effect of predictive mappings on convergence.}
\label{f2_prediction}
\end{minipage}
\begin{minipage}{.48\textwidth}
  \centering
\includegraphics[width=7.5cm]{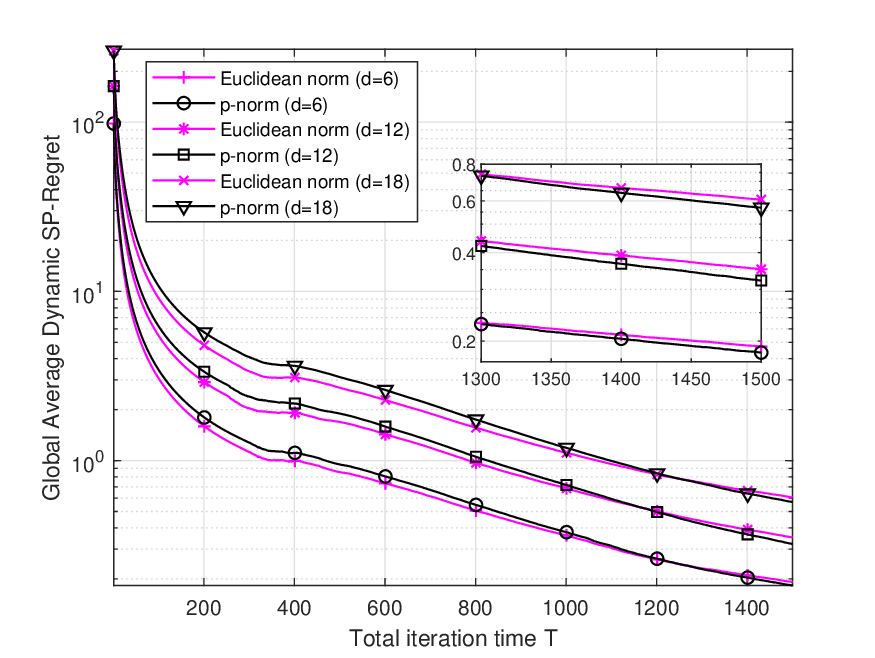}
\caption{Comparison under Euclidean and p-norm.}
\label{f4_2normVSpnorm}
\end{minipage}
\vspace{-0.2cm}
\end{figure}

Secondly, the effect of predictive mappings on convergence of Algorithm \ref{algorithm 1} is analyzed. To avoid noise interference on this result, we set $D_x$ and $D_y$ as 0 and run Algorithm \ref{algorithm 1} with no prediction and with $B_t=P_t, C_t=I_d$, respectively. Fig. \ref{f2_prediction} displays that the convergence performance of the latter algorithm is better than the one without predictions, which implies that the appropriate predictive mappings can effectively enhance the performance of the developed algorithm.

\begin{figure}[!t]
  \centering
\begin{minipage}{.48\textwidth}
  \centering
\includegraphics[width=7.5cm]{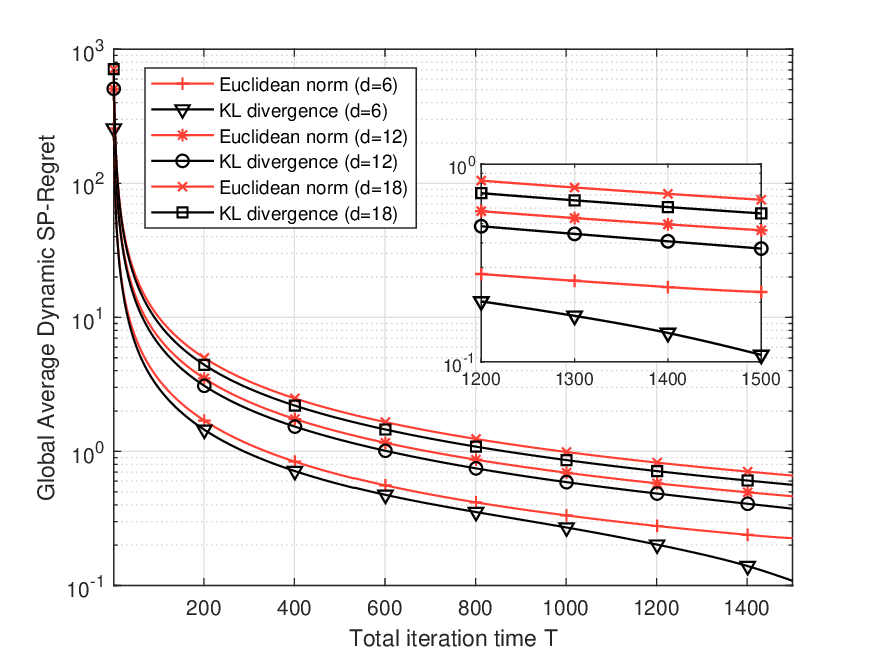}
\caption{The effect of KL divergence on simplex set.}
\label{f4_Sim_2normVsKL}
\end{minipage}
  \centering
\begin{minipage}{.48\textwidth}
  \centering
\includegraphics[width=7.5cm]{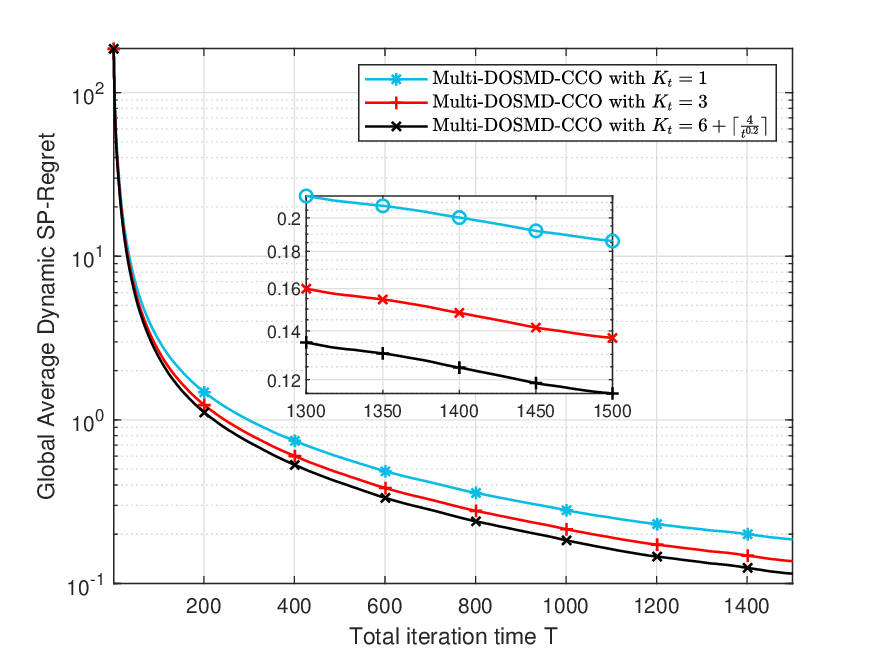}
\caption{The convergence results under different $K_t$.}
\label{f5_multiple}
\end{minipage}
\vspace{-0.2cm}
\end{figure}

Now, we explore the effect of Bregman divergence on the convergence of Algorithm \ref{algorithm 1} in three different dimensions. On the one hand,
the convergence of Algorithm \ref{algorithm 1} using Euclidean norm and $p$-norm is studied and their global ADSPR plots are displayed in Fig. \ref{f4_2normVSpnorm} under the parameters $p=1.85$ and $\lambda_1=0.3$. The results reveal that in the early running stage, the algorithm using Euclidean norm has better convergence than the algorithm using $p$-norm in dimensions $d=6,12,18$,  while as $T$ moves larger, the latter gradually exhibits performance advantages, especially after $T=1200$.
On the other hand, we experiment the convergence of Algorithm \ref{algorithm 1} using KL divergence and Euclidean norm on a simplex constraint. Change $\bsX$ from the current setting to the simplex $\triangle_d$ and keep the setting of $\bsY$. It can be known from \cite{yuan2020distributed} that under $\Psi_{\mathcal{R}}^x (\bsx, \boldsymbol{z})=\sum_{s=1}^d [\bsx]_s \ln (\frac{[\bsx]_s}{[\boldsymbol{z}]_s})$, the mirror step for variable $\bsx_t$ in Algorithm 1 can be written as the explicit solution:
\begin{flalign} \label{mirror step KL}
[\tilde{\bsx}_{i,t}]_s= &\frac{ [\bsx_{i,t}]_s \exp\left(-\alpha_t [\widetilde{\nabla}_{i,t}^x]_s\right)  }{ \sum_{j=1}^d [\bsx_{i,t}]_j \exp\left(-\alpha_t [\widetilde{\nabla}_{i,t}^x]_j\right)  }, s\in [d].
\end{flalign}
Fig. \ref{f4_Sim_2normVsKL} shows that Algorithm \ref{algorithm 1} using KL divergence not only has better convergence performance than the one using Euclidean norm in dimensions $d=6,12,18$, but also is computationally efficient due to its explicit solution in (\ref{mirror step KL}).

\begin{figure}[!t]
  \centering
\begin{minipage}{.48\textwidth}
  \centering
\includegraphics[width=7.5cm]{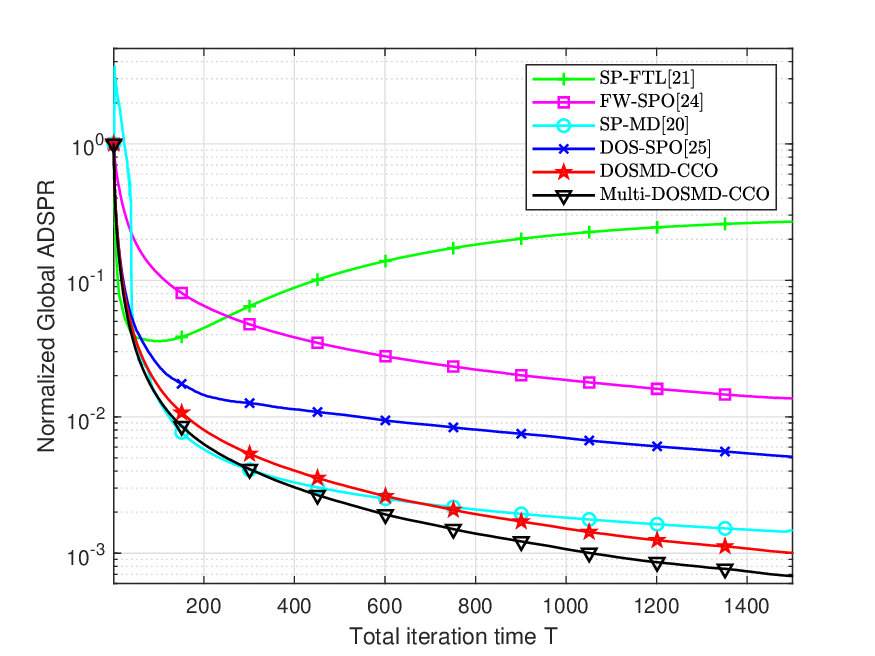}
\caption{Comparisons among different algorithms.}
\label{f6_comparison_Centr_TAC}
\end{minipage}
\vspace{-0.2cm}
\end{figure}

Next, we employ three multiple iteration parameters $K_t=1, 3, 6+\lceil 4/t^{0.2}\rceil$ to investigate the convergence performance of Algorithm Multi-DOSMD-CCO under $n=30$ and analyze the effect of $K_t$ on it.  From Fig. \ref{f5_multiple}, we obtain that 1) When $K_t= 3$ or $ 6+\lceil 4/t^{0.2}\rceil$, the performance of the algorithm are better than the one with $K_t=1$, i.e., Algorithm DOSMD-CCO. In other words, the technique of multiple consensus iterations effectively improves the convergence performance of the primary algorithm; 2) As the parameters $K_t$ become significantly larger, their convergence performance will become better within a certain range, which validates the theoretical results in Corollary \ref{corollary 2}.
Finally, we carry out a rigorous comparison between the proposed algorithms in this paper and distributed online subgradient saddle point optimization (DOS-SPO) algorithm in \cite{10239326}, saddle point mirror descent (SP-MD) algorithm in \cite{ho2019exploiting}, saddle point follow the leader (SP-FTL) algorithm in \cite{rivera2018online} and Frank-Wolfe saddle point optimization (FW-SPO) algorithm in \cite{roy2019online}. To eliminate the performance differences caused by the initial decisions among the studied algorithms, normalized global ADSPR, i.e., $\frac{\frac{1}{n} \sum_j^n \mathcal{M}_j (T)}{\frac{1}{n} \sum_j^n \mathcal{M}_j (1)}$ $\big( \frac{\mathcal{M}_1(T)}{\mathcal{M}_1(1)}$ for centralized algorithms$\big)$, is employed. With $n=30$, $\alpha_t=3/t^{0.35}, \eta_t=15/t^{0.4}$, and the step size requirements from \cite{10239326, ho2019exploiting, roy2019online}, the proposed algorithms display better convergence performance than DOS-SPO algorithm and three centralized algorithms in Fig. \ref{f6_comparison_Centr_TAC} due to the use of predictive mapping and multiple consensus iterations. In addition, note that the regret of SP-FTL algorithm cannot converge, which may be due to the excessive dependence of its decision on the previous loss function.

\section{Conclusions}\label{Conclusions}
In this paper, we have studied the distributed solution for OCCO in a time-varying network. Based on the time-varying predictive technique,  Algorithm DOSMD-CCO has been proposed and analyzed, in which Bregman divergence as a general distance metric and stochastic gradient have been considered. Under the dynamic regret, we have shown that Algorithm DOSMD-CCO guarantees the sublinear convergence for the general convex-concave loss function, provided that $V_T$ is sublinear. Further, to achieve better convergence performance, we also have explored the related multiple consensus version. The obtained results have shown that the appropriate number of consensus iterations can effectively tighten the regret bound to a certain extent.  Finally, the proposed algorithms have been validated and compared through a simulation example involving a target tracking problem. In the future, the promising directions include investigating the nonconvex loss function case and considering a privacy-preserving mechanism.

\begin{appendices}
\section{ PROOFS OF LEMMAS } \label{appendices A}
\subsection{Proof of Lemma \ref{lemma x's diffience}}
To facilitate the analysis, the error variables and their averages are defined as
$
\bse_{i,t}^{x[y]}={\tilde{\bsx}[\tilde{\bsy}]}_{i,t}-{\bsx[\bsy]}_{i,t},\ \bse_{avg,t}^{x[y]}=\frac{1}{n} \sum_{i=1}^n \boldsymbol{e}_{i,t}^{x[y]}
$
where $\bsx[\bsy]$ represents variable $\bsx$ or $\bsy$ in $[\cdot]$. Let $\Pi_{B[C]}(t,s)=B[C]_t B[C]_{t-1}\cdots B[C]_s$, $\forall t\geq s\geq 1$.

According to Algorithm \ref{algorithm 1}, by using the double stochasticity of $A_t$ and a recursive method, we have that
\begin{flalign}
\bsx[\bsy]_{i,t+1}&= \sumjn [\Phi(t,1)]_{ij} \Pi_{B[C]}(t,1) \bsx[\bsy]_{j,1} \nonumber \\
&+ \sum_{\tau=1}^t \sumjn  [\Phi(t, \tau)]_{ij}  \Pi_{B[C]}(t,\tau) \bse_{j,\tau}^{x[y]},\label{lemma x's diffience proof 3-1}\\
\bsx[\bsy]_{avg,t+1}
&=\Pi_{B[C]}(t,1) \bsx[\bsy]_{avg,1} +\sum_{\tau=1}^t \Pi_{B[C]}(t,\tau)  \bse_{avg,\tau}^{x[y]}. \label{lemma x's diffience proof 3-2}
\end{flalign}

Combining (\ref{lemma x's diffience proof 3-1}) and (\ref{lemma x's diffience proof 3-2}), we achieve
\begin{flalign}\label{lemma x's diffience proof 3}
&\left\| \bsx[\bsy]_{i,t+1} - \bsx[\bsy]_{avg,t+1}\right\|\nonumber \\
 &\leq  \sum_{j=1}^n \left|  [\Phi (t,1)]_{ij}- \frac{1}{n} \right| \| \Pi_{B[C]}(t,1)\| \|\bsx[\bsy]_{j,1}\| \nonumber \\
&\quad + \sum_{\tau=1}^t \sum_{j=1}^n  \left|  [\Phi (t,\tau)]_{ij}- \frac{1}{n} \right| \left\|\Pi_{B[C]}(t,\tau)\right\|  \left\|\bse_{j,\tau}^{x[y]}\right\|\nonumber \\
&\overset{(a)}{\leq} \Gamma \sigma^{t-1} \sum_{j=1}^n  \|\bsx[\bsy]_{j,1}\|+ \Gamma \sum_{\tau=1}^t \sigma^{t-\tau} \sum_{j=1}^n \left\|\bse_{j,\tau}^{x[y]} \right\|
\end{flalign}
where (a) is established by using Lemma \ref{lemma-graph} and the fact
\begin{flalign}
\| \Pi_{B[C]}(t,1)\| &= \| B_t B_{t-1}\cdots B_1\|   \nonumber \\
&\leq \| B[C]_t\|  \|B[C]_{t-1}\| \cdots \|B[C]_1\|  \nonumber \\
&\leq 1
\end{flalign}
from Assumption \ref{assu DRX Y noexpansive}.
Summing  the result over $i\in[n]$ and $t\in[T]$, we obtain for $T\geq2$ that
\begin{flalign}\label{lemma x's diffience proof 8}
&\sumT \sumn\| \bsx[\bsy]_{i,t} - \bsx[\bsy]_{avg,t}\|\leq \sumn\| \bsx[\bsy]_{i,1} - \bsx[\bsy]_{avg,1}\| \nonumber \\
 &+ n\Gamma \sum_{t=1}^{T-1} \sigma^{t-1} \sumjn \left\|\bsx[\bsy]_{j,1}\right\|  + n\Gamma \sum_{t=1}^{T-1} \sum_{\tau=1}^t \sigma^{t-\tau} \sumjn \left\|\bse_{j,\tau}^{x[y]} \right\|.
\end{flalign}

The terms on the RHS of (\ref{lemma x's diffience proof 8}) can be bounded as follows:
\begin{flalign}
 1)\ &n\Gamma \sum_{t=1}^{T-1} \sigma^{t-1} \sum_{j=1}^n \|\bsx[\bsy]_{j,1}\|  \leq \frac{ n \Gamma}{1-\sigma} \sum_{j=1}^n \|\bsx[\bsy]_{j,1}\|, \\
  2)\ &n\Gamma \sum_{t=1}^{T-1} \sum_{\tau=1}^t \sigma^{t-\tau} \sum_{j=1}^n \left\|\bse_{j,\tau}^{x[y]} \right\|\nonumber \\
 &\leq n \Gamma \left( \sum_{\tau=1}^{T-1} \sigma^{\tau-1} \right) \sum_{t=1}^{T-1} \sum_{j=1}^n \left\|\bse_{j,\tau}^{x[y]} \right\|\nonumber \\
&\leq \frac{n \Gamma}{ 1-\sigma} \sum_{t=1}^{T-1} \sum_{i=1}^n \| \tilde{\bsx}[ \tilde{\bsy}]_{i,t}-\bsx[\bsy]_{i,t} \|\nonumber \\
&\overset{(b)}{\leq} \frac{n \Gamma }{ \varrho_{x[y]} (1-\sigma)} \sum_{t=1}^{T-1}\sum_{i=1}^n \alpha[\eta]_t  \left\|  \widetilde{\nabla}_{i,t}^{x[y]} \right\|_*
\end{flalign}
where (b) is obtained by combining Lemma \ref{lemma xx yy}.

Finally, Lemma \ref{lemma x's diffience} can be established by substituting the obtained bounds into (\ref{lemma x's diffience proof 8}), and combining the fact that $ \mathbb{E}[\|\widetilde{\nabla}_{i,t}^{x[y]} \|_*]  \leq L_{X[Y]}$.
\hfill$\square$
\subsection{Proof of Lemma \ref{lemma x's diffience2}}
Analogous with the proof of Lemma \ref{lemma x's diffience}, it can be obtained that
\begin{flalign}
1)\ &\bsx[\bsy]_{i,t+1}= \sum_{j=1}^n [\Phi^K(t,1)]_{ij} \Pi_{B[C]}(t,1) \bsx[\bsy]_{j,1}\nonumber\\
  &\quad \quad \quad \quad\quad + \sum_{\tau=1}^t \sum_{j=1}^n [\Phi^K (t, \tau)]_{ij}  \Pi_{B[C]}(t,\tau) \bse_{j,\tau}^{x[y]}, \\
  2)\ &\bsx[\bsy]_{avg,t+1}\!=
\Pi_{B[C]}(t,1) \bsx[\bsy]_{avg,1} \!+\!\sum_{\tau=1}^t \Pi_{B[C]}(t,\tau)  \bse_{avg,\tau}^{x[y]} \label{proof lem x's diff2 eq2}
\end{flalign}
where (\ref{proof lem x's diff2 eq2}) follows double stochasticity of  $(A_t)^{K_t}$.
Combining the above equalities, (\ref{phi-K-b2}), and the facts that $
\| \Pi_{B[C]}(t,1)\|\leq 1,
\sigma_1^{\sum_{p=s}^tK_p-1} \leq \sigma_1^{(t-s+1)\underline{K}-1}$, we get
\begin{flalign}
&\| \bsx[\bsy]_{i,t+1}
 - \bsx[\bsy]_{avg,t+1}\| \nonumber\\
 &\leq \Gamma_1 \sigma_1^{t\underline{K}-1} \sum_{j=1}^n \|\bsx[\bsy]_{j,1}\| \!+\! \Gamma_1 \sum_{\tau=1}^t \sigma_1^{(t-\tau+1)\underline{K}-1} \sum_{j=1}^n \left\|\bse_{j,\tau}^{x[y]} \right\|.
\end{flalign}

Summing  this result over $i\in[n]$ and $t\in[T]$, and using $ \mathbb{E}[\|\widetilde{\nabla}_{i,t}^{x[y]} \|_*]  \leq L_{X[Y]}$ and the methods similar to (\ref{lemma x's diffience proof 8}), Lemma \ref{lemma x's diffience2} can be obtained.
\hfill$\square$
\end{appendices}
\balance

%

\bibliographystyle{IEEEtran}
\bibliography{ArXiv_V7TSP_DOSMD_SPP}

\begin{thebibliography}{10}
\providecommand{\url}[1]{#1}
\csname url@samestyle\endcsname
\providecommand{\newblock}{\relax}
\providecommand{\bibinfo}[2]{#2}
\providecommand{\BIBentrySTDinterwordspacing}{\spaceskip=0pt\relax}
\providecommand{\BIBentryALTinterwordstretchfactor}{4}
\providecommand{\BIBentryALTinterwordspacing}{\spaceskip=\fontdimen2\font plus
\BIBentryALTinterwordstretchfactor\fontdimen3\font minus
  \fontdimen4\font\relax}
\providecommand{\BIBforeignlanguage}[2]{{%
\expandafter\ifx\csname l@#1\endcsname\relax
\typeout{** WARNING: IEEEtran.bst: No hyphenation pattern has been}%
\typeout{** loaded for the language `#1'. Using the pattern for}%
\typeout{** the default language instead.}%
\else
\language=\csname l@#1\endcsname
\fi
#2}}
\providecommand{\BIBdecl}{\relax}
\BIBdecl

\bibitem{shalev2011online}
S.~Shalev-Shwartz \emph{et~al.}, ``Online learning and online convex
  optimization,'' \emph{Found. Trends Mach. Learn.}, vol.~4, no.~2, pp.
  107--194, 2011.

\bibitem{hazan2016introduction}
E.~Hazan \emph{et~al.}, ``Introduction to online convex optimization,''
  \emph{Found. Trends Optim.}, vol.~2, no. 3-4, pp. 157--325, 2016.

\bibitem{yuan2024multiO}
D.~Yuan, A.~Proutiere, G.~Shi \emph{et~al.}, ``Multi-agent online
  optimization,'' \emph{Found. Trends Optim.}, vol.~7, no. 2-3, pp. 81--263,
  2024.

\bibitem{li2023survey}
X.~Li, L.~Xie, and N.~Li, ``A survey on distributed online optimization and
  online games,'' \emph{Annu. Rev. Control}, vol.~56, p. 100904, 2023.

\bibitem{zinkevich2003online}
M.~Zinkevich, ``Online convex programming and generalized infinitesimal
  gradient ascent,'' in \emph{Proc. 20th Int. Conf. Mach. Learn.}, 2003, pp.
  928--936.

\bibitem{cao2019online}
X.~Cao and K.~R. Liu, ``Online convex optimization with time-varying
  constraints and bandit feedback,'' \emph{IEEE Trans. Autom. Control},
  vol.~64, no.~7, pp. 2665--2680, 2019.

\bibitem{9462561}
Y.~Zhao, S.~Qiu, K.~Li, L.~Luo, J.~Yin, and J.~Liu, ``Proximal online gradient
  is optimum for dynamic regret: A general lower bound,'' \emph{IEEE Trans.
  Neural Netw. Learn. Syst.}, vol.~33, no.~12, pp. 7755--7764, 2022.

\bibitem{yuan2020distributed}
D.~Yuan, Y.~Hong, D.~W. Ho, and S.~Xu, ``Distributed mirror descent for online
  composite optimization,'' \emph{IEEE Trans. Autom. Control}, vol.~66, no.~2,
  pp. 714--729, 2020.

\bibitem{yi2020distributed}
X.~Yi, X.~Li, T.~Yang, L.~Xie, T.~Chai, and K.~H. Johansson, ``Distributed
  bandit online convex optimization with time-varying coupled inequality
  constraints,'' \emph{IEEE Trans. Autom. Control}, vol.~66, no.~10, pp.
  4620--4635, 2021.

\bibitem{xiong2022distributed}
Y.~Xiong, X.~Li, K.~You, and L.~Wu, ``Distributed online optimization in
  time-varying unbalanced networks without explicit subgradients,'' \emph{IEEE
  Trans. Signal Proc.}, vol.~70, pp. 4047--4060, 2022.

\bibitem{CaoX2022DisCon}
X.~Cao and T.~Ba\c{s}ar, ``Distributed constrained online convex optimization
  over multiple access fading channels,'' \emph{IEEE Trans. Signal Proc.},
  vol.~70, pp. 3468--3483, 2022.

\bibitem{NazariP2022DAdam}
P.~Nazari, D.~A. Tarzanagh, and G.~Michailidis, ``Dadam: A consensus-based
  distributed adaptive gradient method for online optimization,'' \emph{IEEE
  Trans. Signal Proc.}, vol.~70, pp. 6065--6079, 2022.

\bibitem{9239886}
R.~Dixit, A.~S. Bedi, and K.~Rajawat, ``Online learning over dynamic graphs via
  distributed proximal gradient algorithm,'' \emph{IEEE Trans. Autom. Control},
  vol.~66, no.~11, pp. 5065--5079, 2021.

\bibitem{shahrampour2017distributed}
S.~Shahrampour and A.~Jadbabaie, ``Distributed online optimization in dynamic
  environments using mirror descent,'' \emph{IEEE Trans. Autom. Control},
  vol.~63, no.~3, pp. 714--725, 2017.

\bibitem{beznosikov2020distributed}
A.~Beznosikov, V.~Samokhin, and A.~Gasnikov, ``Distributed saddle-point
  problems: Lower bounds, near-optimal and robust algorithms,'' \emph{Optim.
  Methods Softw.}, 2025, {DOI}: 10.1080/10556788.2025.2463986.

\bibitem{kovalev2022accelerated}
D.~Kovalev, A.~Gasnikov, and P.~Richt{\'a}rik, ``Accelerated primal-dual
  gradient method for smooth and convex-concave saddle-point problems with
  bilinear coupling,'' in \emph{Proc. Adv. Neural Inf. Process. Syst.}, 2022,
  pp. {217}25--{217}37.

\bibitem{akimoto2021saddle}
Y.~Akimoto, Y.~Miyauchi, and A.~Maki, ``Saddle point optimization with
  approximate minimization oracle and its application to robust berthing
  control,'' \emph{ACM Trans. Evol. Learn. Optim.}, vol.~2, no.~1, pp. 1--32,
  2022.

\bibitem{ChenJ2012convergenceSPP}
J.~Chen and V.~K.~N. Lau, ``Convergence analysis of saddle point problems in
  time varying wireless systems--{C}ontrol theoretical approach,'' vol.~60,
  no.~1, pp. 443--452, 2012.

\bibitem{zhou2019adaptive}
H.~Zhou, X.~Zeng, and Y.~Hong, ``Adaptive exact penalty design for constrained
  distributed optimization,'' \emph{IEEE Trans. Autom. Control}, vol.~64,
  no.~11, pp. 4661--4667, 2019.

\bibitem{ho2019exploiting}
N.~Ho-Nguyen and F.~K{\i}l{\i}n{\c{c}}-Karzan, ``Exploiting problem structure
  in optimization under uncertainty via online convex optimization,''
  \emph{Math. Program.}, vol. 177, no.~1, pp. 113--147, 2019.

\bibitem{rivera2018online}
A.~Rivera, H.~Wang, and H.~Xu, ``The online saddle point problem and online
  convex optimization with knapsacks,'' \emph{Math. Oper. Res.}, vol.~50,
  no.~1, pp. 1--39, 2025.

\bibitem{wood2023online}
K.~R. Wood and E.~Dall'Anese, ``Online saddle point tracking with
  decision-dependent data,'' in \emph{Proc. 5th Conf. Learn. Dyn. Control,},
  2023, pp. 1416--1428.

\bibitem{cardoso2019competing}
A.~R. Cardoso, J.~Abernethy, H.~Wang, and H.~Xu, ``Competing against {N}ash
  equilibria in adversarially changing zero-sum games,'' in \emph{Proc. 36th
  Int. Conf. Mach. Learn.}, 2019, pp. 921--930.

\bibitem{roy2019online}
A.~Roy, Y.~Chen, K.~Balasubramanian, and P.~Mohapatra, ``Online and bandit
  algorithms for nonstationary stochastic saddle-point optimization,''
  \emph{arXiv preprint arXiv:1912.01698}, 2019.

\bibitem{10239326}
W.~Zhang, Y.~Shi, B.~Zhang, D.~Yuan, and S.~Xu, ``Dynamic regret of distributed
  online saddle point problem,'' \emph{IEEE Trans. Autom. Control}, vol.~69,
  no.~4, pp. 2522--2529, 2024.

\bibitem{besbes2015non}
O.~Besbes, Y.~Gur, and A.~Zeevi, ``Non-stationary stochastic optimization,''
  \emph{Math. Program.}, vol.~63, no.~5, pp. 1227--1244, 2015.

\bibitem{zhang2019distributed}
Y.~Zhang, R.~J. Ravier, M.~M. Zavlanos, and V.~Tarokh, ``A distributed online
  convex optimization algorithm with improved dynamic regret,'' in \emph{Proc.
  IEEE Conf. Decis. Control}, 2019, pp. 2449--2454.

\bibitem{xu2019online}
Y.~Xu, Y.~Jiang, X.~Xie, and D.~Li, ``An online saddle point optimization
  algorithm with regularization,'' in \emph{Proc. IOP Conf. Series:Mater. Sci.
  Eng.}, 2019, p. 052035.

\bibitem{yang2019survey}
T.~Yang, X.~Yi, J.~Wu, Y.~Yuan, D.~Wu, Z.~Meng, Y.~Hong, H.~Wang, Z.~Lin, and
  K.~H. Johansson, ``A survey of distributed optimization,'' \emph{Annu. Rev.
  Control}, vol.~47, pp. 278--305, 2019.

\bibitem{nedic2018distributed}
A.~Nedi{\'c} and J.~Liu, ``Distributed optimization for control,'' \emph{Annu.
  Rev. Control Robot. Auton. Syst.}, vol.~1, pp. 77--103, 2018.

\bibitem{nedic2008distributed}
A.~Nedic, A.~Olshevsky, A.~Ozdaglar, and J.~N. Tsitsiklis, ``Distributed
  subgradient methods and quantization effects,'' in \emph{Proc. IEEE Conf.
  Decis. Control}, 2008, pp. 4177--4184.

\bibitem{LiuC2024distributed}
C.~Liu, K.~H. Johansson, and Y.~Shi, ``Distributed empirical risk minimization
  with differential privacy,'' \emph{Automatica}, vol. 162, 2024, art. no.
  111514.

\bibitem{NiuD2025adual}
D.~Niu, Y.~Hong, and E.~Song, ``A dual inexact nonsmooth newton method for
  distributed optimization,'' \emph{IEEE Trans. Signal Proc.}, vol.~73, pp.
  188--203, 2025.

\bibitem{YangZ2025Differentially}
Z.~Yang, W.~He, and S.~Yang, ``Differentially private distributed optimization
  over time-varying unbalanced networks with linear convergence rates,''
  \emph{IEEE Trans. Signal Proc.}, vol.~73, pp. 1138--1152, 2025.

\bibitem{xu2021distributed}
J.~Xu, Y.~Tian, Y.~Sun, and G.~Scutari, ``Distributed algorithms for composite
  optimization: Unified framework and convergence analysis,'' \emph{IEEE Trans.
  Signal Proc.}, vol.~69, pp. 3555--3570, 2021.

\bibitem{zhang2020distributed}
J.~Zhang, K.~You, and K.~Cai, ``Distributed dual gradient tracking for resource
  allocation in unbalanced networks,'' \emph{IEEE Trans. Signal Proc.},
  vol.~68, pp. 2186--2198, 2020.

\bibitem{mateos2016distributed}
D.~Mateos-N{\'u}nez and J.~Cort{\'e}s, ``Distributed saddle-point subgradient
  algorithms with {Laplacian} averaging,'' \emph{IEEE Trans. Autom. Control},
  vol.~62, no.~6, pp. 2720--2735, 2016.

\bibitem{rogozin2021decentralized}
A.~Rogozin, A.~Beznosikov, D.~Dvinskikh, D.~Kovalev, P.~Dvurechensky, and
  A.~Gasnikov, ``Decentralized distributed optimization for saddle point
  problems,'' \emph{arXiv preprint arXiv:2102.07758}, 2021.

\bibitem{beznosikov2021distributed}
A.~Beznosikov, G.~Scutari, A.~Rogozin, and A.~Gasnikov, ``Distributed
  saddle-point problems under data similarity,'' in \emph{Proc. Adv. Neural
  Inf. Process. Syst.}, 2021, pp. 8172--8184.

\bibitem{qureshi2023distributed}
M.~I. Qureshi and U.~A. Khan, ``A distributed stochastic first-order method for
  strongly concave-convex saddle point problems,'' in \emph{Proc. IEEE Conf.
  Decis. Control}, 2023, pp. 4170--4175.

\bibitem{polyak1987introduction}
B.~T. Polyak, \emph{Introduction to optimization}.\hskip 1em plus 0.5em minus
  0.4em\relax New York, NY, USA: Optimization Software, Inc, 1987.

\bibitem{eshraghi2020distributed}
N.~Eshraghi and B.~Liang, ``Distributed online optimization over a
  heterogeneous network with any-batch mirror descent,'' in \emph{Proc. 37th
  Int. Conf. Mach. Learn.}, 2020, pp. 2933--2942.

\bibitem{jakovetic2014fast}
D.~Jakoveti{\'c}, J.~Xavier, and J.~M. Moura, ``Fast distributed gradient
  methods,'' \emph{IEEE Trans. Autom. Control}, vol.~59, no.~5, pp. 1131--1146,
  2014.

\bibitem{nedic2014stochastic}
A.~Nedic and S.~Lee, ``On stochastic subgradient mirror-descent algorithm with
  weighted averaging,'' \emph{SIAM J. Optim.}, vol.~24, no.~1, pp. 84--107,
  2014.

\bibitem{xiong2022event}
M.~Xiong, B.~Zhang, D.~W.~C. Ho, D.~Yuan, and S.~Xu, ``Event-triggered
  distributed stochastic mirror descent for convex optimization,'' \emph{IEEE
  Trans. Neural Netw. Learn. Syst.}, vol.~34, no.~9, pp. 6480--6491, 2023.

\bibitem{banerjee2005clustering}
A.~Banerjee, S.~Merugu, I.~S. Dhillon, J.~Ghosh, and J.~Lafferty, ``Clustering
  with bregman divergences.'' \emph{J. Mach. Learn Res.}, vol.~6, no.~10, 2005.

\bibitem{bauschke2001joint}
H.~H. Bauschke and J.~M. Borwein, ``Joint and separate convexity of the bregman
  distance,'' \emph{Stud. Comput. Math.}, vol.~8, pp. 23--36, 2001.

\bibitem{dhillon2008matrix}
I.~S. Dhillon and J.~A. Tropp, ``Matrix nearness problems with bregman
  divergences,'' \emph{SIAM J. Matrix Anal. Appl.}, vol.~29, no.~4, pp.
  1120--1146, 2008.

\bibitem{yi2020distributed2}
X.~Yi, X.~Li, L.~Xie, and K.~H. Johansson, ``Distributed online convex
  optimization with time-varying coupled inequality constraints,'' \emph{IEEE
  Trans. Signal Proc.}, vol.~68, pp. 731--746, 2020.

\bibitem{li2021distributed}
J.~Li, C.~Li, W.~Yu, X.~Zhu, and X.~Yu, ``Distributed online bandit learning in
  dynamic environments over unbalanced digraphs,'' \emph{IEEE Trans. Netw. Sci.
  Eng.}, vol.~8, no.~4, pp. 3034--3047, 2021.

\bibitem{zhang2018adaptive}
L.~Zhang, S.~Lu, and Z.-H. Zhou, ``Adaptive online learning in dynamic
  environments,'' in \emph{Proc. Adv. Neural Inf. Process. Syst.}, 2018, pp.
  1330--1340.

\bibitem{ben2009robust}
A.~Ben-Tal, L.~El~Ghaoui, and A.~Nemirovski, \emph{Robust optimization}.\hskip
  1em plus 0.5em minus 0.4em\relax Princeton, NJ, USA: Princeton Univ. Press,
  2009.

\bibitem{zhang2023communication}
S.~Zhang, S.~Choudhury, S.~U. Stich, and N.~Loizou, ``Communication-efficient
  gradient descent-accent methods for distributed variational inequalities:
  Unified analysis and local updates,'' \emph{arXiv preprint arXiv:2306.05100},
  2023.

\end{thebibliography}

\end{document}